\documentclass[10pt,a4paper]{amsart}

\usepackage{amsthm}
\usepackage{latexsym}
\usepackage[mathscr]{eucal}
\usepackage{amscd}
\usepackage{amssymb}
\usepackage{amsmath}
\usepackage{amsfonts}
\usepackage[mathscr]{eucal}
\usepackage[sans]{dsfont}
\usepackage{stmaryrd}
\usepackage[all]{xy}
\usepackage{color}
\def\1{{\mathds 1}}

\newcommand\su[2]{\overset{#2}{\underset{#1}\sum}\,}
\newcommand\osu[2]{\overset{#2}{\underset{#1}\oplus}\,}
\newcommand\bosu[2]{\overset{#2}{\underset{#1}\bigoplus}\,}
\newcommand\cu[2]{\overset{#2}{\underset{#1}\cup}\,}
\newcommand\ca[2]{\overset{#2}{\underset{#1}\cap}\,}
\newcommand\pru[2]{\overset{#2}{\underset{#1}\prod}\,}

\DeclareMathOperator{\Hom}{Hom}
\DeclareMathOperator{\Tor}{Tor}
\DeclareMathOperator{\Mod}{Mod}
\DeclareMathOperator{\mo}{mod}
\DeclareMathOperator{\Ker}{Ker}
\DeclareMathOperator{\Ext}{Ext}
\DeclareMathOperator{\End}{End}
\DeclareMathOperator{\Aut}{Aut}
\DeclareMathOperator{\Specm}{Specm}
\DeclareMathOperator{\Ob}{Ob}
\DeclareMathOperator{\gl}{gl}
\DeclareMathOperator{\ssl}{sl}
\DeclareMathOperator{\id}{id}

\DeclareMathOperator{\gkdim}{GKdim}
\DeclareMathOperator{\supp}{supp}
\DeclareMathOperator{\Spec}{Spec}
\DeclareMathOperator{\imm}{Im}
\DeclareMathOperator\St{St}
\renewcommand\Im{\imm}

\newcommand\bC{\mathbb{C}}
\newcommand\bD{\mathbb{D}}
\newcommand\bH{\mathbb{H}}
\newcommand\bF{\mathbb{F}}
\newcommand\bZ{\mathbb{Z}}

\newcommand\cA{\mathscr{A}}
\newcommand\cB{\mathscr{B}}
\newcommand\cC{\mathscr{C}}
\newcommand\cI{\mathscr{I}}
\newcommand\cK{\mathscr{K}}
\newcommand\cL{\mathscr{L}}
\newcommand\cM{\mathscr{M}}
\newcommand\cO{\mathscr{O}}
\newcommand\cX{\mathscr{X}}
\newcommand\kX{\mathcal{X}}
\newcommand\kY{\mathcal{Y}}
\newcommand\kN{\mathcal{N}}

\renewcommand\sl{{\mathsf l}}
\newcommand\sm{\mathsf {m}}
\newcommand\sn{\mathsf {n}}
\renewcommand\sp{\mathsf{p}}
\newcommand\sq{\mathsf {q}}

\newcommand\rM{\mathrm{M}}

\renewcommand{\k}{\Bbbk}
\newcommand\Ga{\Gamma}
\newcommand\La{\Lambda}
\newcommand\ga{\gamma}

\newcommand\Om{\Omega}

\newcommand\vi{\varphi}

\renewcommand\r[2]{#2_{#1}}
\renewcommand\l[2]{{}_{#1}#2}
\newcommand\lr[3]{{}_{#1}#3_{#2}}
\newcommand\Lr[3]{\,{}^{#1}#3_{#2}}
\newcommand\lR[3]{{}_{#1}\,#3^{#2}\,}

\newcommand\ds{\displaystyle}
\newcommand\myto{\longrightarrow}
\newcommand{\xar}[1]
{{\,{\xrightarrow{\ \ #1\ \ }\,}}}

\newtheorem{definition}{Definition}
\newtheorem{theorem}{Theorem}[section]
\newtheorem{corollary}{Corollary}[section]
\newtheorem{lemma}{Lemma}[section]
\newtheorem{proposition}{Proposition}[section]
\newtheorem{remark}{Remark}[section]

\newcommand\incirc[1]{{\xy*+{#1}*\cir{}\endxy}}

\newcounter{numberofremark}
\newcounter{numberofremarkA}
\begin{document}

\author{Vyacheslav Futorny }
\author{Serge Ovsienko  }
\title[Fibers of characters
\  ]{Fibers of characters in Gelfand-Tsetlin categories }
\address{Instituto de Matem\'atica e Estat\'\i stica,
Universidade de S\~ao Paulo, Caixa
Postal 66281, S\~ao Paulo, CEP
05315-970, Brasil}
\email{futorny@ime.usp.br}
\address{
Faculty of Mechanics and
Mathematics, Kiev Taras Shevchenko
University, Vla\-di\-mir\-skaya 64,
00133, Kiev, Ukraine}
\email{ovsienko@zeos.net }


\begin{abstract}
We solve the problem of extension of characters of commutative
subalgebras in associative (noncommutative) algebras for a class
of subrings (Galois orders) in skew group rings. These results can
be viewed as a noncommutative analogue of liftings of prime ideals
in the case of integral extensions of commutative rings. The
proposed approach can be applied to the representation theory  of
many infinite dimensional algebras including  universal enveloping
algebras of reductive Lie algebras,  Yangians and finite
$W$-algebras. In particular, we develop a
theory of Gelfand-Tsetlin modules for $\gl_n$.
Besides classification results we characterize their categories in
the generic case extending the classical results on $\gl_2$.
\end{abstract}

\maketitle

\tableofcontents

\label{will-be}

\section{Introduction}
\label{section-introduction}
The functors of restriction  onto subalgebras and  induction from
subalgebras are important tools in the representation theory. The effectiveness of these tools depends upon the choice of a subalgebra. Denote by $\Specm A$ ($\Spec A$) the space of maximal (prime)
ideals in $A$, endowed with the Zarisky topology.
In the classical commutative algebra setup,  an integral extension
$A\subset B$ of two commutative rings   induces a  map
$\vi:\Spec B\rightarrow\Spec A$, whose fibers are
non-empty for every point of $\Spec A$ (e.g. $A=B^G$, where $G$ is a finite subgroup of the
automorphism group of $B$). In particular, every character of
$A$ can be extended to a
character of  its integral extension $B$. Moreover, if $B$ is finite over $A$
then all fibers $\vi^{-1}(I)$, $I\in \Spec A$ are finite, and hence, the number of extensions of a character of $A$ is finite.   The
Hilbert-Noether theorem provides an example of such situation with
$B$ being the symmetric algebra on a finite-dimensional vector
space $V$ and $A=B^G$, where $G$ is a
finite subgroup of $GL(V)$.

The primary goal of this paper is to generalize these results to
 "semi-commu\-ta\-ti\-ve" pairs $\Ga\subset U$ where $U$ is an
associative (noncommutative) algebra over a base field $\k$ and
$\Ga$ is  an integral domain.
 The canonical embedding $\Ga\subset U$
induces a functor from the category of $U$-modules which are
direct sums of finite dimensional $\Ga$-modules
(\emph{Gelfand-Tsetlin modules with respect to} $\Ga$) to the
category of torsion $\Ga$-modules. This functor induces a
``multivalued function'' from  $\Specm \Ga$ associating to an
ideal $ \sm\in\Specm \Ga$    the \textit{fiber} $\Phi(\sm)$  of
left maximal ideals of $U$ that contain $\sm$. Our goal is to find
natural sufficient conditions for the fibers  to be non-empty and
finite for any point in $\Specm \Ga$. On the other hand, for a
maximal left ideal $I\subset U$ such that $U/I$ is a
Gelfand-Tsetlin module, it is interesting to investigate its
support in $\Specm \Ga$ (i.e. the set of $\sm\in\Specm \Ga$, such
that $\Ga/\sm$ is a subfactor of $U/I$ as a $\Ga$-module) and find
the multiplicity of $\Ga/\sm$ in $U/I$.

A motivation for the study of such pairs $(U, \Ga)$  comes from
the representation theory.    The classical framework of
\textit{Harish-Chandra modules}  (\cite{d}, Ch. 9) is related to a
pair of a reductive Lie algebra $\mathcal F$ and its reductive
subalgebra $\mathcal F'$, where $U$ and $\Gamma$ are their
universal enveloping algebras respectively. A more general concept of
Harish-Chandra modules (related to a pair $(U,\Gamma)$) was introduced in
\cite{dfo:hc}.

 The case when $U$ is the universal enveloping
algebra of a reductive finite dimensional Lie algebra and $\Ga$ is
the universal enveloping algebra of a Cartan subalgebra leads to
the theory of Harish-Chandra modules with respect to this Cartan
subalgebra, commonly known as generalized weight  modules.
Classification of such simple  modules is well known for $\gl_2$ and for any simple finite-dimensional
Lie algebra for
  modules with finite-dimensional
weight spaces, due to Fernando \cite{Fe} and
Mathieu \cite{Ma}.
 It remains an open problem in general. To
approach this classification problem, the full subcategory of
weight Gelfand-Tsetlin  $U(\gl_n)$-modules with respect to the
\emph{Gelfand-Tsetlin subalgebra}  was introduced in
\cite{dfo:gz}. This class  is based on natural properties of a
Gelfand-Tsetlin basis for finite-dimensional representations of
simple classical Lie algebras \cite{g-ts}, \cite{zh:cg},
\cite{m:gtsb}.

Gelfand-Tsetlin subalgebras were considered in
 \cite{FM} in  connection with the solutions of the
Euler equation, in \cite{Vi} in  connection with subalgebras of
maximal Gelfand-Kirillov dimension in the universal enveloping
algebra of a simple Lie algebra, in \cite{KW-1}, \cite{KW-2} in
connection with classical mechanics, and also in \cite{Gr} in
connection with general hypergeometric functions on the Lie group
$GL(n, \bC)$. A similar approach was used by Okunkov and
Vershik in their study of the representations of the symmetric group
$S_n$ \cite{OV}, with $U$ being the group algebra of $S_n$ and $\Ga$
being the maximal commutative subalgebra generated by the
Jucys-Murphy elements
$$(1i)+\ldots+ (i-1 i)\,, \quad \,i=1,\ldots, n.$$  In this case the
elements of $\Specm \Ga$ parametrize  basis of irreducible
representations of $U$.
Recent advances in the
representation theory of Yangians (\cite{FMO}) and  finite $W$-algebras (\cite{FMO1}) are also
based on similar techniques.

What is the intrinsic reason
of the existence of Gelfand-Tsetlin formulae and of the
successful study of Gelfand-Tsetlin representations of
various classes of algebras? This question led us to the introduction in
\cite{fo-GaI} of concepts
of  \emph{Galois rings} and  \emph{ Galois orders}  in
invariant skew (semi)group rings.

\emph{For the rest of the paper we assume that  $\Ga$ is a commutative domain, $K$  the field of fractions of $\Ga$, $K\subset L$ a finite Galois extension,   $\cM\subset
\Aut L$ a submonoid closed under conjugation by the elements of the Galois group $G=G(L/K)$. We will always assume that if for any $m_{1},m_{2}\in \cM$ from
$$m_{1}|_{K}=m_{2}|_{K}$$
follows $m_{1}=m_{2}$. Such monoid $\cM$ we call \emph{separating} with respect to $K$.}

The group $G$ acts on the skew semigroup ring $L*\cM$ via the action $(lm)^g=l^g m^g$, where $m^g=g^{-1}m g$. By $\cK $ we denote the subring  of $G$-invariants   $(L*\cM)^{G}\subset L*\cM.$

\begin{definition}
\label{definition-definition-Galois-ring-and-order}
A \emph{Galois ring} $U$ over $\Gamma$ is a finitely generated over $\Gamma$  subring in $\cK$ such that $KU=UK=\cK$. A Galois ring $U$  over   $\Ga$ is called \emph{right (respectively left) order}  if for any finite dimensional right (respectively left) $K$-subspace $W\subset \cK$ (respectively $W\subset
\cK$), $W\cap U$ is a finitely generated right (respectively
left) $\Ga$-module. A Galois ring is \emph{order}{} if it is both right and left order  (\cite{fo-GaI}).
\end{definition}
Galois orders are natural versions of ``noncommutative orders'' in skew semigroup rings of invariants.
In  comparison with  the classical notion of an order  we note, that $\Ga\subset U$ is not central, \emph{but maximal commutative subalgebra}.

A class of Galois orders  includes in particular the following
subrings in the corresponding skew group rings (\cite{fo-GaI}):
 Generalized Weyl algebras over integral domains with infinite order automorphisms, e.g. $n$-th Weyl algebra
$A_n$, quantum plane, $q$-deformed Heisenberg algebra, quantized
Weyl algebras, Witten-Wo\-ro\-no\-wicz algebra among the others
(\cite{ba}, \cite{bavo}); the universal enveloping algebra of
$\gl_n$ over the Gelfand-Tsetlin subalgebra (\cite{dfo:gz},
\cite{dfo:hc}),  associated shifted Yangians and finite
$W$-algebras (\cite{FMO}, \cite{FMO1});   certain rings of
invariant differential operators on torus. In Section~
\ref{section-Preliminaries} we present some necessary  facts about
Galois orders.

In this paper we develop a representation theory of Galois orders.
The main tool for our investigation of  categories of
representation is a technique from \cite{dfo:hc}. In Section~
\ref{section-Gelfand-Tsetlin-categories} we give a detailed
exposition  of the results from  \cite{dfo:hc} adapted for
 case of a commutative subalgebra considered in this paper.

Last two sections are devoted to the representation theory of
Galois orders. We emphasize  that the theory of Galois orders unifies the
representation theories of universal enveloping algebras and
generalized Weyl algebras. Our main result establishes  sufficient
conditions for  the fiber $\Phi(\sm)$ to be nontrivial and finite.
Let $\ell_{\sm}$ be any lifting of $\sm$ to the integral closure
of $\Ga$ in $L$, $\cM_{\sm}$ the stabilizer of $\ell_{\sm}$ in
$\cM$. Note that the group $\cM_{\sm}$ is defined uniquely up to $G$-conjugation, hence its cardinality is correctly defined.

Our main result is the following
\medskip

\noindent{\bf Theorem A.} \textit{Let $\Gamma$ be a commutative domain
which is finitely generated as a $\k$-algebra,  $U$ a
Galois ring over $\Ga$, $\sm\in
\Specm \Ga$. Suppose that $\cM_{\sm}$ is finite.
\begin{itemize}
\item\label{fiber-nontrivial} If $U$ is a right
  Galois order  over
$\Ga$ then the fiber $\Phi(\sm)$ is non-empty.
\item\label{fiber-finite} If $U$ is a
  Galois order  over
$\Ga$ then the fiber $\Phi(\sm)$ is finite.
\end{itemize}}

\medskip

The methods we use in the proof of these results allow us to show
that Gelfand-Tsetlin modules for a large  class of  Galois orders
over $\Ga$  have similar ``nice'' properties. For any $\sm\in
\Specm \Ga$ with finite $\cM_{\sm}$ we obtain an effective
estimate for the number of isomorphism classes of simple
Gelfand-Tsetlin modules $M$ whose support contains $\sm$ and for
the dimension of generalized weight spaces $M(\sm)$  (Theorem~
\ref{theorem-finite-dimensionity-of-component}). In particular,
for $U=U(\gl_{n})$ these numbers are limited by $2!\ldots (n-1)!$.

In Section~
\ref{section-Representations-of-general-position} we collected the
results related to Gelfand-Tsetlin representations in ``general
position''. The corresponding blocks in module category have a
unique simple representation, and such a  block is equivalent to
the category of finite dimensional representations of a completed
commutative algebra. As an application  we obtain a version of the
Harish-Chandra theorem for Galois orders, (Theorem B,
\eqref{enumerate-massive-faithfull-action}).

\medskip

\noindent{\bf Theorem B.}\emph{ Let $\Gamma$  a commutative domain which is finitely
generated and normal as a $\k$-algebra, $U$ a Galois order over
$\Ga$ with $\cM$ being a group.
\begin{enumerate}
\item\label{enumerate-massive-simple}
There exists a massive subset $W\subset \Specm \Ga$ such
that for any $\sm\in W$, $|\Phi(\sm)|=1$ and hence there exists a
unique simple $U$-module $L_{\sm}$ whose support contains $\sm$.
\item\label{ext-category} The extension category generated by $L_{\sm}$ contains all
indecomposable modules whose support contains $\sm$ and it is
equivalent to the category of modules over the completion ${\Ga}_{\sm}$.
\item\label{enumerate-massive-faithfull-action} If $\cM$ is a group, then for any nonzero $u\in U$ there exists a massive set of
non-isomorphic simple Gelfand-Tsetlin $U$-modules on which $u$ acts
nontrivially.
\end{enumerate}
}

Note, that in the case $U=U(\ssl_{2})$ the structure of the
category of Gelfand-Tsetlin modules is well known (see, e.g.
\cite{dr}). This category splits into a direct sum of blocks, and
each block is equivalent to the category of finite dimensional
representations of complete algebras over ${\k[c]}_{\sm},$ where
$c\in U$ is the Casimir element and $\sm\in\Specm \k[c]$ is a
maximal ideal which acts on the corresponding block nilpotently.
These algebras are presented by the quivers
%{\ar@{-->}@/^1.4pc/|<<<<{\alpha_{22}} (30,22)*+{}; (29,-6)*+{} }
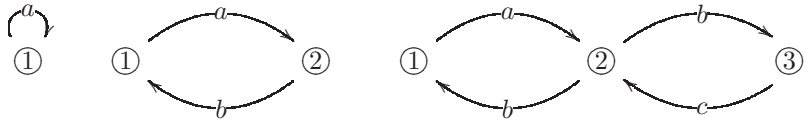
\begin{equation*}
%\label{equation-sl-2-gz-quivers}
\xymatrix@!=400pc{
\incirc{1}\ar@(ul,ru)[]|{\ds a}
}\quad\quad
\xymatrix@!=1pc{
\incirc{1}\ar@{->}@/^{1.5pc}/|{\ds a}[rr]&&\incirc{2}\ar@{->}@/^{1.5pc}/|{\ds b}[ll]
} \quad\quad
\xymatrix{
\incirc{1}\ar@{->}@/^{1.5pc}/|{\ds a}[rr]&&\incirc{2}\ar@{->}@/^{1.5pc}/|{\ds b}[ll] \ar@{->}@/^{1.5pc}/|{\ds b}[rr]&&\incirc{3}\ar@{->}@/^{1.5pc}/|{\ds c}[ll].
\quad
 }
\end{equation*}
The first quiver corresponds to the case of generic blocks, while
the second case corresponds to the blocks which contain Verma
module, but does not contain a finite dimensional module. The
third quiver corresponds to the blocks, containing a finite
dimensional $\ssl_{2}$-module. Besides,  in  this last case holds
the relation $ab=cd.$ \footnote{The problem of classification of
finite dimensional modules over this algebra is a famous ``Gelfand
problem'', \cite{ge}.}  In the situation of an arbitrary Galois
order (even in the case of the universal enveloping algebra of $\ssl_{n}$) we are far from such
detailed description of the blocks of the category of
Gelfand-Tsetlin modules.  Nevertheless we are able to characterize
the structure of generic blocks.  Theorem  \textbf{B} above
describes the generic blocks of the category of Gelfand-Tsetlin
modules, the situation is analogous to the case of $\ssl_{2},$
namely generic blocks are equivalent to module categories for the algebra of formal power series.

A step in the direction of a characterization of the structure of
the blocks of the category of Gelfand-Tsetlin modules is
Proposition \ref{proposition-break}.  In
Section~\ref{section-repres-integral} we establish sufficient
conditions under which a block of a Gelfand-Tsetlin category
contains a finite number of non-isomorphic simple modules
(Corollary~\ref{corollary-general-estimate-of-length}).
Note that an analogue of this statement (see \cite{dfo:hc},
Theorem 32) allows to prove the statement, that the "subgeneric"
blocks of Gelfand-Tsetlin modules for $U(\gl_{n})$ are described
by a finite quiver with relations. The structure of the blocks,
and closely related question of the finiteness of length of the
left $U-$module $U/U\sm$ we will discuss in a subsequent paper.

We note here an important connection which arose in the case when
$U=U(\gl_{n})$  and $\Ga\subset U$ is a Gelfand-Tsetlin
subalgebra. In this case an important role is played by the
variety of so-called strongly nilpotent matrices (\cite{O}). In
 was shown in \cite{O} that this variety is a complete intersection.
 In particular, it implies that $U$ is free both as a right and as a left
$\Ga$-module (\cite{fo3}). Kostant and Wallach (\cite{KW-1},
\cite{KW-2})  introduced a generalization of the variety of
strongly nilpotent matrices and revealed a deep relation between
this variety and the hamiltonian mechanics.  A connection between
 Gelfand-Tsetlin representations of $U$ and  the structure of the
Kostant-Wallach variety is evidently important and should be a
topic of a further study.

In this paper we apply the  theory only to Lie algebras
of type $A$, but we believe it can be extended to other types.
This technique was used in \cite{FMO1} to address the
classification problem of irreducible Gelfand-Tsetlin modules for
finite \emph{$W$-algebras} and \emph{shifted Yangians} associated
with $\gl_n$ and to prove an analogue of the Gelfand-Kirillov
conjecture for these algebras.

\section{Preliminaries}
\label{section-Preliminaries}

\subsection{Notations}
\label{subsection-List of notations}

All fields in the paper contain the base  field $\k$, which
is algebraically closed of characteristic $0$. All algebras
in the paper are $\k$-algebras.  If $A$ is an associative ring then by $A-\mo$ we denote the
category of finitely generated left $A$-modules. Let $\cC$
be a $\k$-category, i.e. all $\Hom_\cC$-sets are endowed
with a structure of a $\k$-vector space and all the
compositions maps are $\k$-bilinear. The category of
$\cC$-modules $\cC-\Mod$ is defined as the category of
$\k$-linear functors $M:\cC\myto \k-\Mod$, where $\k-\Mod$
is the category of $\k$-vector spaces.
The category of locally finitely generated $\cC$-modules we denote by $\cC-\mo$. Let $G$ be a group, $X$ a $G$-set, then by $X/G$ we denote the corresponding factor-set and by $X^{G}$ the set of $G$-invariants. For a set $X$ by $|X|$ we denote the cardinality of $X.$ Let $G$ be a group, $H\subset G$ a subgroup. Then \emph{the notation $\su{g\in G/H}{}F(g)$ will mean, that element $g$ runs a set of representatives of $G/H$ under the assumption that the sum does non depend on the choice of these representatives.}

\subsection{Integral extensions}
Let $A$ be an integral commutative domain, $K$
its field of fractions and
$\tilde{A}$ the integral closure of
$A$ in $K$. The ring $A$ is called \emph{normal}
if $A=\tilde{A}$. The following is standard

\begin{proposition}\label{basics-integral}
Let $A$ be a  normal noetherian
ring, $K\subset L$ a finite Galois
extension, $\bar{A}$ is the
integral closure of $A$ in $L$.
Then $\bar{A}$ is a finite
$A$-module.
\end{proposition}

Let $\imath:A\hookrightarrow B$ be an integral extension. Then it
induces a surjective map $\Specm B\rightarrow \Specm A$ ($\Spec
B\myto \Spec A$). In particular, for any character
$\chi:A\rightarrow \k$ there exists a character
$\tilde{\chi}:B\rightarrow \k$ such that $\tilde{\chi}|_A=\chi$.
If, in addition, $B$ is finite over $A$, i.e. finitely generated
as an $A$-module, then the number of different characters of $B$
which correspond to the same character of $A$, is finite.

\begin{corollary}\label{cor-finite-extension}  (\cite{s}, Ch. III, Prop. 11, Prop. 16)
If $A$ is a finitely generated
$\k$-algebra then for any character
$\chi:A\rightarrow \k$ there exists
finitely many characters
$\tilde{\chi}:\bar{A}\rightarrow
\k$ such that
$\tilde{\chi}|_A=\chi$.
\end{corollary}

The following statement is probably well known but in \cite{fo-GaI} we
include the proof for the convenience of the reader. Note that we consider Proposition \ref{prop-nonsingular-commutative-lifting} as a motivation for introducing the notion of a Galois order.

\begin{proposition} (\cite{fo-GaI})
\label{prop-nonsingular-commutative-lifting} Let
$i:A\hookrightarrow B$ be an embedding of integral domains
with a regular $A$. Assume the induced morphism of
varieties $i^{*}:\Specm B\rightarrow \Specm A$ is
surjective (e.g. $A\subset B$ is an integral extension). If
$b\in B$ and $a b\in A$ for some nonzero $a\in A$ then
$b\in A$.
\end{proposition}

\subsection{Skew (semi)group rings}
\label{subsection-Skew-group-rings-and-skew-group-categories}
Let $R$ be a ring, $\cM$ a submonoid of $\Aut R$.
The \emph{skew semigroup ring}, $R*\cM$,  is a free left
$R$-module with a basis $\cM$ and with the
multiplication
$$(r_{1}m_{1})\cdot (r_{2}m_{2})= (r_{1}r_{2}^{m_{1}}) (m_{1}
m_{2}),\quad m_{1},m_{2}\in \cM,\ r_{1},r_{2}\in R.$$

If $x\in R*\cM$ and $m\in \cM$ then denote by $x_{m}$ the element
in $R$ such that $x=\su{m\in \cM}{} x_{m} m$. Set
$$\supp x=\{m\in\cM|x_{m}\ne0\}.$$
If a finite group $G$ acts by automorphisms on $R$ and by
conjugations on $\cM$ then $G$ acts on $R*\cM$.  Denote by
 $(R*\cM)^{G}$  the invariants under this action.
 Then $x\in (R*\cM)^{G}$ if and only if
$x_{m^{g}}=x^{g}_{m}$ for $m\in\cM,g\in G$.

For $\vi\in\Aut R$ and $a\in R$ set $H_{\vi}=\{h\in
G|\vi^{h}=\vi\}$ and
\begin{equation}
\label{equation-definition-a-phi} [a\vi]:=\sum_{g\in
G/H_{\vi}}a^{g}\vi^{g}\in (R*\cM)^{G}.\footnote{Recall, that the notation  $g\in
G/H_{\vi}$ means that $g$ runs over a set of representatives of cosets from $G/H_{\vi}$ and the result does not depend on a choice of these representatives.}
\end{equation}

Then

\begin{equation}
\label{equation-decomp-LM-G-in-invariant}
\begin{split}
&(R*\cM)^{G}=\bigoplus_{\vi\in \cM/G}(R*\cM)^G_{\vi},\
\text{ where } \\ &(R*\cM)^G_{\vi}= \big\{[a\vi]\,|\,a\in
R^{H_{\vi}}\,\big\}.
\end{split}
\end{equation}

We will use the following formulae:
\begin{align} \label{space-structure}
&\gamma \cdot [a
\vi]=[(a\gamma) \vi],\ [a \vi]\cdot \gamma= [(a\gamma^{\vi}) \vi],
\ \ga\in R^{G},\\ \nonumber
&\text{ and } [a\vi]=
[\vi a^{\vi^{-1}}].
\end{align}

\begin{remark}
\label{remark-squarebracket-gamma-generated-bimodule}
The formulae \eqref{space-structure} means, that as left $R$-module $R[a\vi]R$  is isomorphic to
$RR^{\vi}$ and as a right $R$-module  $R[a\vi]R$ is isomorphic to $RR^{\vi^{-1}}$.
\end{remark}

%For $a,b\in R^{H_{\vi}}$, $\gamma\in R^G$ denote

%\begin{align}
%\label{align-def-a-vi-b} &[a\vi
%b]=\sum_{g\in
%G/H_{\vi}}a^{g}\vi^{g}b^{g},\\
%\nonumber & \ga [a\vi b]=[(\ga a)\vi b]=[a\vi (b\ga^{\vi^{-1}})],\
%[a\vi b]\ga=[(\ga^{\vi} a)\vi b]=[a\vi (b\ga)].
%\end{align}

\subsection{Galois rings}
\label{subsection-Galois-rings}
The notations and results Subsection \ref{subsection-Skew-group-rings-and-skew-group-categories} we will use in the case when  $R=L$ is a field, $K\subset L$ is a finite Galois extension of fields, $G=G(L/K)$ its Galois group. We will denote by $\imath$  the canonical embedding $K\hookrightarrow L.$
Recall, that we will use the notations introduced before Definition \ref{definition-definition-Galois-ring-and-order}.

%
%A \emph{Galois ring} $U$ over $\Gamma$ is a finitely generated over $\Gamma$
%subring of the ring  $\cK=(L*\cM)^{G}$ of invariants such that
%$KU=UK=\cK$.
%

% Let $\ell_{\sm}$ be any lifting of $\sm$ to the integral closure of $\Ga$ in $L$ , $\cM_{\sm}$ the stabilizer of $\ell_{\sm}$ in $\cM$. Since $\sm$ defines  $\cM_{\ell_{\sm}}$ up to conjugacy,
%by abuse of notations we denote it by $\cM_{\sm}$.
%
%
%
%
%
%Denote  $K=L^{G}$.
%
%\begin{remark}
%\label{remark-all-attributes}
%Since in the text below we will work with a fixed Galois ring,  we will use the following agreement. We assume, that by definition every Galois $U$ ring is endowed with the following data, presented in the definition of a Galois ring:
%$$ \Ga\subset K, K\subset L, G=G(L/K), \cM\subset \Aut L \text{ is a separating monoid } $$
%$$ g\cM g^{-1}=\cM, g\in G; \cK=(L*\cM)^{G},  $$
%and with a canonical embedding $i:U\hookrightarrow \cK$, identical on $\Ga$.
%
%
%$K,L,G,\cM, \cK, V(S)$, to explain the notation $|S/G|$ - all ingredients!
%
%$I_{L}(\Ga)$
%
%\end{remark}
%
%
%

Let $S\subset \cM$ be a finite $G$-invariant subset, $V\subset \cK$ a $\Ga$-subbimodule. Then we introduce   the following $\Ga$-subbimodule in $V$
\begin{equation}
\label{equation-definition-of-support-subsets}
V(S)=\{x\in V\mid \supp x \subset S\}
\end{equation}
We need the following simple fact.

\begin{lemma}
\label{lemma-dimension-with-fixed-support} In the assumption above the right and the left dimensions of $\cK(S)$ coincides with $|S/G|$.
\end{lemma}

In particular, for any $\vi\in\cM$ holds
\begin{equation}
\label{equation-dim-K} \dim_{K}\cK_{\vi}=
[L^{H_{\vi}}:K]=[G:H_{\vi}]=|G{\vi}|.
\end{equation}
$\cK_{\vi}$ is irreducible as a $K$-bimodule. Such bimodule is denoted in \cite{fo-GaI} by $V(\vi).$

Note also, that since $\Ga$ acts torsion-free on $\cK$ both from the left and from the right, we obtain that the canonical maps
\begin{align}
\label{align-cK-obtained-by-tensoring}
&e_{r}:U\otimes_{\Ga}K\myto \cK,\ u\otimes x\mapsto ux, u\in U, x\in K;\\ \nonumber
&e_{l}:K\otimes_{\Ga}U\myto \cK,\ x\otimes u\mapsto xu, u\in U, x\in K
\end{align}
isomorphism of $\Ga-K$  and $K-\Ga$-bimodules correspondingly.

Recall, that the monoid $\cM\subset \Aut L$ from definition of a Galois ring we assume to be  \emph{separating} (with respect to $K$), i.e. if for any $m_{1},m_{2}\in \cM$ the equality $m_{1}|_{K}=m_{2}|_{K}$ implies $m_{1}=m_{2}$. An automorphism $\vi:L\myto L$ is called     \emph{separating} (with respect to $K$) if the monoid generated by $\{\vi^{g}\,|\,g\in G\}$ in $\Aut L$ is separating.

\begin{lemma}\label{lemma-separ-monoid}(\cite{fo-GaI}, Proposition 2.3) Let  $\cM$ be a separating monoid with respect to
    $K$. Then
\begin{enumerate}
\item  \label{enum-def-sep-action-nonunit} $\cM\cap G=\{e\}.$
\item\label{enum-def-of-sep-action-act} For any
    $m\in\cM,m\ne e$ there exists $\ga\in K$ such that
    $\ga^{m}\ne\ga$.

\item\label{enum-def-of-sep-action-galois} If $G
    m_{1}G=G m_{2} G$ for some $m_{1},m_{2}\in\cM$,
    then there exists $g\in G$ such that
    $m_{1}=m_{2}^{g}$.
\item\label{enum-def-of-sep-action-group-case} If $\cM$ is a group, then the statements \eqref{enum-def-sep-action-nonunit}, \eqref{enum-def-of-sep-action-act}, \eqref{enum-def-of-sep-action-galois} are equivalent and each of them implies that $\cM$ is separating.
\end{enumerate}
\end{lemma}

Let $U$ be a \emph{Galois $\Ga$-ring}{}
over $\Ga$ (cf. Introduction). Denote by $i$ the canonical embedding $i:\Ga\hookrightarrow U$.

\begin{lemma}\label{l81}(compare \cite{fo-GaI}, Lemma 4.1)
\label{lemma-Ka-viK-and-right-left-basis}
Let $u\in U$ be nonzero element, $$T=\supp u, \,u=\ds\sum_{m\in
T} [a_m m].$$ Then
$$K(\Ga u \Ga)=(\Ga u \Ga)K= K uK=\bigoplus_{m\in T} V(a_{m}m),$$
where $V(a_{m}m)=K[a_m m]K$ is an irreducible $K$-bimodule.
\end{lemma}

In particular it shows that for every $m\in \cM$, $U$
contains the elements $[b_{1}m],$ $\dots,$ $[b_{k}m]$ where
$b_{1},\dots,b_{k}$ is a $K$-basis in $L^{H_{m}}$. Let $e\in \cM$ be the unit element and $U_{e}=
U\cap Le$.

\begin{theorem}(\cite{fo-GaI}, Theorem 4.1)
\label{theorem-shat-algebras} Let $U$ be a Galois ring over $\Gamma$. Then
\begin{enumerate}

\item \label{enum-GA-is-G-invariant} $U_{e}\subset K$.

\item \label{enum-shat-is-maximal} $U\cap K$ is a maximal
commutative $\k$-subalgebra in $U$.

\item \label{enum-shat-center} The center $Z(U)$ of  $U$
equals $U\cap K^{\cM}$.
\end{enumerate}
\end{theorem}

\subsection{Galois orders and Harish-Chandra subalgebras}
\label{subsection-Integral-Galois-algebra}
In this section we recall
the basic properties of    Galois orders following
\cite{fo-GaI}. For simplicity we only consider right
Galois orders. Let $M$ be a right $\Gamma$-submodule of a Galois
order $U$ over $\Gamma$. Set
$$\bD_{r}(M)=\{u\in U\,| \text{ there exists } \gamma\in \Gamma,
\gamma\ne 0 \text{ such that } u
\cdot\gamma\in M\}.$$

We have the following characterization of    Galois
orders.

\begin{proposition}(\cite{fo-GaI}, Proposition 5.1)
\label{prop-integral-Galois-algebra} A Galois ring $U$  over  a
noetherian $\Ga$ is a \emph{right order} if and only if for every
finitely generated right $\Gamma$-module $M\subset U$, the right
$\Gamma$-module $\bD_{r}(M)$ is finitely generated.
\end{proposition}

In particular, if $U$ is right integral then $\Ga\subset U_e$ is an
integral extension and $U_e$ is a normal ring.

Recall that  ${\Gamma}$ is called a \emph{Harish-Chandra subalgebra} in $U$ if $\Ga u\Ga$ is finitely generated both as a
left and as a right $\Ga$-module for any $u\in U$  \cite{dfo:hc}. We will also say that $\Gamma$ is a
\emph{right (left) Harish-Chandra subalgebra} if $\Ga u\Ga$ is finitely generated  as a
right (left) $\Ga$-module for any $u\in U$. Note, that these property is enough to check for some set of generators of the ring $U$ over $\Ga.$ Harish-Chandra subalgebras of Galois rings have the following properties.

\begin{proposition}\label{integral-hchandra}
If $\, U$ is a right (left)  Galois order over a noetherian
$\Ga$ then for any $m\in \cM$ holds $m^{-1}(\Ga)\subset \bar{\Ga}$
($m(\Ga)\subset \bar{\Ga}$).
\end{proposition}

\begin{proposition}
\label{proposition-order-and-integrality-Harish-Chandra}
Assume  $\Ga$ is
finitely generated algebra over $\k$, $U$ is a Galois
ring. Then $\Ga$ is a Harish-Chandra subalgebra in $U$ if and only
if $m\cdot \bar{\Ga}= \bar{\Ga}$ for every $m\in\cM$.
\end{proposition}

\begin{proposition}(\cite{fo-GaI}, Corollary 5.3)
\label{proposition-order-implies-Harish-Chandra}
If $U$ is a  Galois order  over  $\Ga$ and
 $\Ga$ is a
noetherian $\Bbbk$-algebra
 then $\Ga$ is a Harish-Chandra
subalgebra in $U$.
\end{proposition}

The next Lemma is a main technical tool in our investigation of representations of Galois orders. Let $S\subset \cM$ be a finite $G$-invariant subset. Denote
$$U(S)=\{u\in U\, |\,\supp u\subset S\}.$$
For every $f\in \Gamma$ consider $f_{S}^r\subset
\Gamma\otimes_{\k}K$  as follows
\begin{equation}
\label{equation-relations-in-as-in-gz}
f_{S}^r=\prod_{s\in S}(f\otimes 1 -
1\otimes
f^{s^{-1}})=\sum_{i=0}^{|S|}
f^{|S|-i}\otimes T_{i}, \
(T_{0}=1).
\end{equation}

Similarly we define  $f_{S}^l=\pru{s\in S}{}(f^s\otimes
1 - 1\otimes f)\in K\otimes_{\k}\Ga.$

\begin{lemma}(\cite{fo-GaI}, Lemma 5.2)
\label{lemma-relation-as-in-gz-acts-zero}  If $m^{-1}(\Ga)\subset \bar{\Ga}$ ($m(\Ga)\subset \bar{\Ga}$ respectively) for all
$m\in \cM,$ then for any $G$-invariant subset $X\subset \cM$ and $f_X=f_X^r$ ($f_X=f_X^l$ respectively) holds $f_{X}\in \Ga\otimes_{\k}\Ga.$
Besides for a $G$-invariant subset $S\subset \cM$ holds the following.

\begin{enumerate}

\item\label{enum-f-S-kills-U-S}  $u\in U(S)$ if and
    only if $f_{S}\cdot u=0$ for every $f\in\Gamma$.

\item \label{enum-f-complement-to-S-throws-U-S} If
    $T=\supp u \setminus S$ then $f_{T}\cdot u\in
    U(S)$ for every $f\in\Gamma$.

\item\label{enum-f-S-mult-on-standard} If
    $f_{S}=\ds\sum_{i=1}^{n}f_{i}\otimes g_{i}$,
    $[am]\in L*\cM$ then $$f_{S}\cdot
    [am]=[(\ds\sum_{i=1}^{n} f_{i}g_{i}^{m}a)m] = [\prod_{s\in S }(f -
    f^{ms^{-1}})am].$$

\item\label{enum-f-S-mult-components-homomorphism} Let
    $S$ be  a $G$-orbit and $T$ an $G$-invariant subset
    in $\cM$. The $\Ga$-bimodule homomorphism
    $P^{T}_{S}(=P^{T}_{S}(f)):U(T)\myto U(S)$,
    $u\mapsto f_{T\setminus S}\cdot u,f\in \Gamma$ is
    either zero or $\Ker P_{S}^{T}=U(T\setminus S)$.

\item\label{enum-f-S-mult-split-components-mono}
Let $S=S_{1}\sqcup\dots\sqcup
S_{n}$ be the decomposition of $S$
in $G$-orbits and
$P_{S_{i}}^{S}:U(S)\myto U(S_{i})$,
$i=1,\dots,n$ are defined in
\eqref{enum-f-S-mult-components-homomorphism}
nonzero homomorphisms.  Then the
homomorphism
\begin{equation}
\label{equation-bimodule-distinguishing-components}
P^{S}:U(S)\myto\bigoplus_{i=1}^{n}
U(S_{i}),\
P^{S}=(P_{S_{1}}^{S},\dots,P_{S_{n}}^{S}),
\end{equation}
is a monomorphism.

\end{enumerate}
\end{lemma}

We have the following equivalent conditions for Galois ring to be a Galois order.

\begin{theorem}(cf. \cite{fo-GaI}, Theorem 5.1)
\label{theorem-integrality-equivalence-of-defs} Let $U$ be a Galois
ring over  a noetherian $\Ga$ and assume that $\Ga$ is a right (left)
 Harish-Chandra $\k$-subalgebra of $U$. Then the following
statements are equivalent.

\begin{enumerate}

\item\label{enum-equiv-of-defs-int-right}
$U$ is right (respectively left)   Galois order.

\item\label{enum-equiv-of-defs-int-right-any-s}
$U(S)$ is finitely generated right
(respectively left) $\Ga$-module
for any finite $G$-invariant
$S\subset \cM$.

\item\label{enum-equiv-of-defs-int-right-elements}
$U(G\cdot m)$ is finitely generated
right (respectively left)
$\Ga$-module for any $m\in\cM$.

\end{enumerate}
\end{theorem}

\begin{theorem}(cf. \cite{fo-GaI}, Theorem 5.2)
\label{theorem-for-involutive-algebra-right-is-left} Let
$U$ be a Galois ring over  a noetherian $\Ga$ and $\cM$  a subgroup of $\Aut L$.
\begin{enumerate}

\item
\label{enum-involutive-right-integral} If $U_{e}$ is integral
extension of $\Gamma$ and $m^{-1}(\Gamma)\subset \bar{\Gamma}$
(respectively $m(\Gamma)\subset \bar{\Gamma}$), then $U$ is a right
(respectively left)   Galois order if and only if $\,U_{e}$ is an integral
extension of $\Gamma$.

\item
\label{enum-involutive-U-0-is-integral}
If $U_{e}$ is an integral extension of $\Gamma$ and $\Gamma$ is a Harish-Chandra $\k$-subalgebra in $U$, then $U$ is  Galois order over $\Ga$.
\end{enumerate}
\end{theorem}

The corollary below gives us a  converse statement to {\bf{Theorem~{B}}}.
\begin{corollary}
\label{corollary-existence-of-lifting-implies-Galois-order}
 Let $U\subset L*\cM$ be a Galois ring over a noetherian $\Ga$.
If every character $\chi:\Gamma\myto \k$ extends to a representation
of $U$ then $U_{e}\subset\bar{\Gamma}\cap K$. If in addition $\cM$
is a group and $\Ga$ is a Harish-Chandra subalgebra then $U$ is
a Galois order.
\end{corollary}

\begin{proof}
If $\chi$ extends to a representation of $U$, then it
extends to a representation of $U_{e}\subset K$ in
particular.  Proposition \ref{prop-nonsingular-commutative-lifting}  implies that $U_{e}$ belongs to the integral
closure of $\Gamma$ in $K$. The second statement follows
immediately from
Theorem~\ref{theorem-for-involutive-algebra-right-is-left}.
\end{proof}

The next corollary sets a bridge between the theory of (noncommutative) Galois orders and commutative case in Proposition~\ref{prop-nonsingular-commutative-lifting}.

\begin{corollary}
\label{corollary-criterion-integra-algebra}(\cite{fo-GaI}, Corollary 5.6) Let $U\subset L*\cM$ be
a Galois ring over noetherian $\Ga$, $\cM$ a group and $\Ga$ a
normal $\k$-algebra. Then the following statements are equivalent
\begin{enumerate}
\item \label{enum-algebra-is-integral} $U$ is a Galois order.

\item \label{enum-divisor-criterion-integrality} $\Ga$ is a Harish-Chandra subalgebra and, if for $u\in U$ there exists a nonzero $\ga\in \Ga$ such that $\ga
    u\in \Ga$ or $u\ga\in \Ga$, then $u\in\Gamma$.
\end{enumerate}
\end{corollary}

%\begin{proposition}
%\label{prop-reduction-to-commutative} Assume $U$ is endowed with increasing exhausting filtration, such that the %associated graded algebra $\gr{U}$ is a commutative domain.  If the canonical
%embedding $\imath: \gr{\Gamma} \hookrightarrow \gr{U}$
%induces an epimorphism $$\imath^{*}:\Specm \gr{U} \to
%\Specm \gr{\Gamma}$$ then $U$ is a Galois order over $\Ga$.
%\end{proposition}

\section{Gelfand-Tsetlin categories}
\label{section-Gelfand-Tsetlin-categories}

\subsection{Motivation}
\label{subsection-Motivations-Gelfand-Tsetlin-categories}

The constructions of this section are the main tools we will use
to investigate the class of Gelfand-Tsetlin $U$-modules. First
such constructions appeared in \cite{dfo:hc} in general setting, but for our purposes
we consider here a special case of a
commutative subalgebra $\Ga$ and present it  in  details.
In this section we just
assume that $\Ga$ is a commutative Harish-Chandra subalgebra in a
finitely generated associative algebra $U$.

Before going into details we give  some motivation of the constructions below.
Let $U$ be a finite dimensional associative algebra over $\k$, $R\subset U$ its Levi subalgebra, $\Ga\subset U$ the center of $R$. Then $\Ga=\osu{i=1}{n}\k e_{i}$, where $\{e_{1},\dots,e_{n}\}$ is the complete (i.e. $e_{1}+\dots+e_{n}=1$) family of mutually orthogonal idempotents. Obviously, $\Ga\subset U$ is a Harish-Chandra subalgebra.

This data allows to present $U$ as a ring of the $n\times n$ matrices of the  form \begin{equation}
\label{equation-pearce-decomposition}\left(e_{j}Ue_{i}\right)_{i,j=1,\dots,n}, \quad  u\longmapsto \left(\begin{array}{ccc}
                     e_{1}ue_{1} & \dots & e_{1}ue_{n} \\
                     \vdots & \ddots & \vdots \\
                     e_{n}ue_{1} & \dots & e_{n}ue_{n}
                   \end{array}
                   \right)
 \end{equation}
 with the standard multiplication. This presentation is called the \textit{two-sided Pierce decomposition of $U$.} Besides, we can associate with $U$ a $\k$-linear category
\begin{equation*}
\cA=\cA(U;\Ga) \text{ where } \Ob \cA=\{1,\dots,n\},  \, \cA(i,j)=e_{j}U e_{i},
\end{equation*}
and the composition of morphism is defined by the multiplication in $U$. One simple but important remark is the existence of the canonical isomorphism
\begin{equation}
\label{equation-module-equivalence-for-finite-dimensional}
U-\Mod\simeq \cA-\Mod.
\end{equation}
If $\Ga=R,$ equivalently, $U$ is a basic (or Morita reduced) algebra, then the category $\cA$  is presented usually as a quiver with relations. This presentation is the key feature in the study of finite  dimensional representations of $U$ (see e.g. \cite{dk}, \cite{gr} for details).

The definition of $\cA$ can be rewritten as follows. As objects of
$\cA$ one can consider $\Specm \Ga=$ $\{\sm_{1},\dots,\sm_{n}\}$,
where $\sm_{i}$ is the kernel of the projection of $\Ga$ onto $\k
e_{i}.$ Besides, set $\cA(\sm_{i},\sm_{j})=\Ga/\sm_{j}\otimes_{\Ga}
U\otimes_{\Ga} \Ga/\sm_{i}$ with the composition of morphisms
induced by the multiplication in $U$. The construction  of the category $\cA$ in \cite{dfo:hc} can be considered as a generalization of the two sided Pierce decomposition. The construction below (see Definition~\ref{equation-hom-in-A}) is a special case  of a commutative Harish-Chandra subalgebra $\Ga\subset U$, where $U$ is not  necessary finite dimensional. As above, we associate with the pair $\Ga\subset U$ a category $\cA$ with  $\Ob \cA=\Specm \Ga$. Unfortunately,  there is no  equivalence of categories of $U$-modules and $\cA$-modules. Instead we have a weaker result for the full subcategory of Gelfand-Tsetlin $U$-modules (see Theorem~\ref{theorem-HC-mod-equiv-A-mod}).

\subsection{Gelfand-Tsetlin modules}
\label{subsection-Gelfand-Tsetlin-modules}
We assume $U$  an algebra over $\k$ and $\Ga\subset U$ is a commutative finitely generated  subalgebra.
The following is the key notion of the paper.
\begin{definition}
\label{definition-of-GZ-modules} A finitely generated $U$-module
$M$ is called \emph{Gelfand-Tsetlin module (with respect to
$\Gamma$)}{} provided that the restriction $M|_{\Gamma}$ is a
direct sum of  $\Gamma$-modules
\begin{equation}
\label{equation-Gelfand-Tsetlin-module-def}
M|_{\Gamma}=\bigoplus_{\sm\in\Specm \Gamma}{}M(\sm),
\end{equation}
where $$M(\sm)=\{v\in M| \sm^{k}v=0 \text{ for some }k\geq 0\}.$$
\end{definition}

When for all $\sm\in \Specm {\Gamma}$ and all $x\in M(\sm)$ holds
$\sm x=0$ such Gelfand-Tsetlin module $M$ is called \emph{weight
module (with respect to $\Gamma$)}{}. More generally, for a left
(right) $\Ga$-module $X$ and $\sm\in\Specm \Ga$ we call an element
$x\in X$ left (right) \textit{$\sm$-nilpotent}, provided that
$\sm^{k}x=0$  (x$\sm^{k}=0$) for some $k\geq1.$

All
Gelfand-Tsetlin modules form a full abelian and closed with
respect to extensions subcategory
\index{$\bH(U,{\Gamma})$} $\bH(U,{\Gamma})$
 in $U-{\mo}$. A full subcategory of $
\bH(U,{\Gamma})$ consisting of weight Gelfand-Tsetlin we denote by
$ \bH W(U,{\Gamma})$.

The \textit{support} of a
Gelfand-Tsetlin module $M$ is a set $$\supp M=\{\sm\in\Specm
{\Gamma}\mid M({\sm})\ne 0\}.$$ For $D\subset$ $\Specm {\Gamma}$ denote by $
\bH(U,{\Gamma},D)$ the full subcategory in $\bH(U,{\Gamma})$ formed by $M$ such that $\supp M\subset$
$D$. For a given $\sm\in \Specm {\Gamma}$ we denote
by $\chi_\sm:\Gamma\myto \Gamma/\sm$ the
corresponding character of $\Gamma$. Conversely, for a character $\chi:\Ga\myto \k$ denote $\sm_{\chi}=\Ker \chi$, so  we will identify the set of
all characters of $\Gamma$ with $\Specm \Ga$. If
there exists a simple Gelfand-Tsetlin module $M$ with
$M(\sm)\ne 0$ then we say that the character $\chi_\sm$
\emph{lifts}  to $M$.

\emph{From now on we  assume that $\Ga$ is a Harish-Chandra subalgebra in $U$.} For a $\Ga$-bimodule $V$ any pair $(\sm,\sn)\in \Specm {\Gamma} \times \Specm {\Gamma}$ and $m,n\geq0$
 we will use the following notations:
\begin{align*}
& \l {\sn^{n}}{V}=V/\sn^{n}V, \quad \r {\sm^{m}}{V}=V/V\sm^{m},\quad \lr{\sn^{n}}{\sm^{m}}{V}=V/(\sn^{n}V+V\sm^{m}).
\end{align*}

For $a\in U$ and $V={\Gamma} a {\Gamma}$ denote the first two bimodules above by $L_{a,m}$ and $R_{a,n}$ respectively.

\begin{lemma}
\label{lemma-corollaries-of-hc-conditions}
In the assumption above holds the following.

\begin{enumerate}
\label{enumerate-simple-HC-properties} \item
\label{item-finite-Ga-length} The modules $L_{a,m},R_{a,n}$  are
finite dimensional.

\item\label{item-finite-Ga-input-output-in-lattice} Any from the
following three conditions
\begin{enumerate}
\label{enumerate-three-def}
\item $\Gamma/\sn $ is the subquotient of $L_{a,1}$ as a left $\Ga$-module.

\item $\Gamma/\sm $ is the subquotient of $R_{a,1}$ as a right $\Ga$-module.

\item $\Ga/ \sn \otimes_{\Ga} \Ga a\Ga \otimes_{\Ga} \Ga/\sm \ne0.$
\end{enumerate}
defines the same set  $X_{a}$ of  pairs $ (\sm,\sn)\in \Specm
{\Gamma} \times \Specm {\Gamma}.$

\item \label{item-finite-Ga-input-output-in-ideal} For $\sm,\sn\in \Specm \Ga$ define the  set $$X_{a}(\sm)=\{\sn\in\Specm
\Ga \mid(\sm,\sn)\in X_{a}\}$$
and the set
 $$X^{a}(\sn) =\{\sm\in\Specm
\Ga \mid(\sm,\sn)\in X_{a}\}.$$ Then both $X_{a}(\sm)$ and
$X^{a}(\sn) $ are finite. Besides, the kernels of simple
subquotients of all $L_{a,m}$ (respectively $R_{a,n}$)
belong to $X_{a}(\sm)$ (respectively $X^{a}(\sn)$).

\item\label{equation-1}   Let $M$ be a Gelfand-Tsetlin $U$-module,
$\sm\in \Specm \Ga$. Then
    \begin{equation}
    a M(\sm) \ \subset
   {\sum_{\sn\in X_a(\sm)}}M(\sn).
    \end{equation}

\item
    \label{equation-2}  Let $M$ be a Gelfand-Tsetlin $U$-module,  $\pi_{\sn}:M\myto M(\sn)$,  $\sn\in \Specm \Ga$
     the canonical projection and $\sm\not\in X^{a}(\sn)$. Then
     \begin{equation}
    \pi_{\sn}(a M(\sm))=0.
 \end{equation}

 \item\label{enumerate-induced-from-f-d-is-GZ}
 If $X$ is a finite-dimensional ${\Gamma}-$module then
$U\otimes_{{\Gamma}} X$ is a Gelfand-Tsetlin module.

\end{enumerate}
\end{lemma}

\begin{proof}

We will prove the statements for $L_{a,m}$. The case of $R_{a,n}$
is analogous.

Since $\Ga$ is finitely generated, $\dim_{k}\Ga/\sm^{m}< \infty$.
Then $L_{a,m}\simeq\Ga a \Ga\otimes_{\Ga} \Ga/\sm^{m}$ is
finite dimensional, since $\Ga a \Ga$ is finitely generated as a
right $\Ga$-module.  This proves \eqref{item-finite-Ga-length}. If $L_{a,\sm}=\osu{k=1}{t} L_{k}$ is a decomposition
into a direct sum of indecomposable left $\Ga$-modules, then for
every $k=1,\dots,t$ there exists $\sn_{k}\in \Specm \Ga$ and
$n_{\k}\geq1$, such that $\sn_{k}^{n_{k}}L_{k}=0$.  In particular,
the subquotients $L_{k}$ are isomorphic to $\Ga/\sn_{k}$. On the
other hand, $\Ga/\sn\otimes_{\Ga} L_{k}\simeq L_{k}/\sn L_{k}$ is
nonzero if and only if $\sn=\sn_{k}$, which, together with
\eqref{item-finite-Ga-length}, proves
\eqref{item-finite-Ga-input-output-in-lattice} and
\eqref{item-finite-Ga-input-output-in-ideal}.
To prove \eqref{equation-1} consider any $x\in M(\sm)$. Then there
exists $m>1$, such that $\sm^{m}x=0$. It follows that the left
$\Ga$-submodule  $\Ga a \Ga x\subset M$ is the factor of
$L_{a,\sm^{m}}$. Then the statement  follows from
\eqref{item-finite-Ga-input-output-in-ideal}. The statement
\eqref{equation-2} is proved analogously.
To show \eqref{enumerate-induced-from-f-d-is-GZ} it is enough to
consider the case $\dim_{\k} V=1.$ But then the statement  follows from
\eqref{item-finite-Ga-length}.
\end{proof}

Denote by $\Delta$
the minimal equivalence on $\Specm {\Gamma}$ containing all
$X_a$, $a\in U$ and by $ {\Delta}(U, {{\Gamma}})$ the set
of the $\Delta-$equivalence classes on $\Specm {\Gamma}$.

\begin{lemma}
\label{lemma-ext-disjointness} Let $X,X'$ be Gelfand-Tsetlin
modules,  $\supp X\subset$ $ D,$ $\supp
X'\subset$ $ D'$, where $D$ and $D'$ are different
classes of $\Delta$-equivalence. Then
$$\Hom_{U}(X,X')=0,\ \Ext_{U}^{1}(X,X')=0.$$
\end{lemma}

\begin{proof}
Obviously, $\Hom_{\Ga}(X,X')=0$ all the moreover $ \Hom_{U}(X,X')=0.$
It is enough to prove that every exact sequence in $U-\mo$
\begin{equation*}
0\longrightarrow X'\xrightarrow{ \sigma } Y \xrightarrow{ \pi } X\longrightarrow 0
\end{equation*}
 splits. Since $D\cap D'=\varnothing$, this sequence splits \textit{uniquely} as a sequence of $\Ga$-modules,
 thus  we can assume $$Y(\sm)=X'(\sm)\oplus X(\sm), \sm\in \supp Y.$$ Since $D\cap D'=\varnothing$, in this sum
 either $X'(\sm)=0$
or $X(\sm)=0.$ For $a\in U$  holds $a X'(\sm)\subset X'$ and $a
X(\sm)\subset X$, by Lemma
\ref{lemma-corollaries-of-hc-conditions}, \eqref{equation-1}.
Hence $X'$ and $X$ are $U$-submodules in $Y$.
\end{proof}

Immediately from Lemma \ref{lemma-corollaries-of-hc-conditions}, \eqref{equation-1}, \eqref{equation-2} and Lemma \ref{lemma-ext-disjointness} we obtain the following
\begin{corollary}
\label{corollary-splitting-gz-modules-in-blocks}
$$\bH(U,{\Gamma})=
\bigoplus_{D\in
\Delta(U,{\Gamma})}\bH(U,{\Gamma},D).$$
\end{corollary}

The subcategory $\bH(U,{\Gamma},D)$, where $D\in
\Delta(U,{\Gamma})$ will be called a \textit{block} of $\bH(U,\Gamma)$.

\subsection{The category $\cA$}
\label{subsection-The-category-ScAS}

For a $\Ga$-bimodule $V$ denote by $_{\sn}\hat{V}_{\sm}$ the  $I$-adic  completion of $\Ga\otimes_{\k}\Ga$-module $V$, where  $I\subset \Ga\otimes \Ga$ is a maximal ideal $I=\sn\otimes\Ga+\Ga\otimes \sm$, in other words
\begin{equation*}
 _{\sn}\hat{V}_{\sm} \ = \
\ds\varprojlim_{n,m} \lr{\sn^n}{ \sm^m}{V}.
\end{equation*}

Let $B$ be a $\Ga$-bimodule, satisfying \emph{the Harish-Chandra
condition: every finite generated bimodule in $B$ is
finitely generated  both from the left and from the right.} Denote
by $F(B)$ the set of finitely generated $\Ga$-subbimodules  in
$B$. Note, that if  $V,W\in F(B)$ and $V\subset U$, then the
canonical embedding induces a monomorphism
$_{\sn}V_{\sm}\hookrightarrow _{\sn}W_{\sm}.$ It allows us to give
the following definition.

The \textit{finitary completion} $_{\sn}\widetilde{B}_{\sm}$  of the bimodule $B$ by the ideal $I=\sn\otimes\Ga+\Ga\otimes \sm$  is a $(\widehat{\Ga\otimes \Ga})_{I}$ module
\begin{equation}
\label{equation-finitary-completed}
_{\sn}\widetilde{B}_{\sm}=\varinjlim_{V\in F(B)} {}_{\sn}\hat{V}_{\sm}.
\end{equation}

For any $V\in F(B),$ $\sm,\sn\in \Specm \Ga,$ $m,n>0$ we have a family of  $\Ga$-bimodule morphisms
\begin{equation*}
_{\sn}\hat{V}_{\sm}\myto  _{\sn^{n}}V_{\sm^{m}}\myto  _{\sn^{n}}B_{\sm^{m}},
\end{equation*}
where the first is the canonical map from the projective limit and the second is induced by the embedding $V\subset B.$  This family defines a homomorphism
\begin{equation*}
\Psi:_{\sn}\widetilde{B}_{\sm}=\varinjlim_{V\in F(B)} {}_{\sn}\hat{V}_{\sm}\myto \varprojlim_{n,m} {}_{\sn^n}B_{ \sm^m} = _{\sn}\hat{B}_{\sm}.
\end{equation*}
If  $B$ is finitely generated as $\Ga$-bimodule then $\Psi$ is an isomorphism.

For $\sm \in$ $\Specm
{\Gamma}$ denote by 
$\displaystyle \hat{\Gamma}_{\sm}=$ $\displaystyle
\lim_{\leftarrow m}{\Gamma}/\sm^m$  the
completion of ${\Gamma}$.

\begin{definition}
\label{definition-of-category-cA}\label{equation-hom-in-A}
Define a category \index{${\cA}$} ${\cA}$ $=$ ${\cA}_{U,{\Gamma}}$ with objects $\Ob{\cA}=$ $\Specm {\Gamma}$. The space of morphisms from
$\sm$ to $\sn$ is defined as the completion of $U$,  i.e.
\begin{equation*}
 {\cA}(\sm,\sn) \ = \ _{\sn} \widetilde{U}_{\sm}\ \big(=
\varinjlim_{V\in F(U)} \varprojlim_{n,m} {} _{\sn^n}V_{ \sm^m}.\big)
\end{equation*}
\end{definition}

The spaces  $\cA(\sm,\sn)$ are endowed with the standard topology
defined by the limits. Besides, any $\cA(\sm,\sn)$ is endowed with
a canonical  structure of a completed
$\hat{\Ga}_{\sn}- \hat{\Ga}_{\sm}$-bimodule. For $\sl\in\Specm\Ga$  set
\begin{align*}
& \Lr{\sl}{\sm^{m}}{V}=\{ \bar{a}
\in \r{\sm^{m}}{V}\mid \sl^{l}\bar a=0 \text{ for some } l\geq 1\},\\ \nonumber
& \lR{\sn^{n}}{\sl}{V}=\{ \bar{b}
\in \l {\sn^{n}}{V}\mid \bar b\sl^{l}=0 \text{ for some } l\geq 1\}.
\end{align*}

For $V\subset U$ consider $\r{\sm^{m}}{V}$ as a left $\Ga$-module.
By Lemma \ref{lemma-corollaries-of-hc-conditions},
\eqref{item-finite-Ga-input-output-in-ideal} for every $a\in V$
there exists a finite set $X_{a}(\sm)=\{\sn_{1},\dots,\sn_{k}\}$ and
$N\geq 1, $ $N=N(a,\sm),$ such that $\sn_{1}^{N}\dots\sn_{k}^{N}$
annihilates the class $\bar{a}$ of $a$ in $V/V\sm^{m}.$ Hence by
the Chinese remainder theorem
\begin{align}\label{align-right-left-root-HC-subalgebra-decomposition} \nonumber
&\r{\sm^{m}}{V}=\bigoplus_{\sl\in X_{\sm}} \Lr{\sl}{\sm^{m}}{V},  \text{ where } X_{\sm}=\bigcup_{a\in V}X_{a}(\sm)\\
&\text{ and }\\ \nonumber
&\l{\sn^{n}}{V}=\bigoplus_{\sl\in X^{\sn}} \lR{\sn^{n}}{\sl}{V},  \text{ where } X^{\sn}=\bigcup_{a\in V}X^{a}(\sn).
\end{align}

\begin{lemma}
\label{lemma-simplest-properties-of-cA}
\begin{enumerate}
\item\label{enumerate-properties-of-cA-disjont-points}
Let $D,$ $D'\subset \Specm \Ga$ be two different classes of $\Delta$-equi\-valence, $\sm\in D$, $\sn\in D'.$ Then $\cA(\sm,\sn)=0.$

\item\label{enumerate-properties-of-cA-disjoint-cA-D}
We have a decomposition of $\cA$ into a direct sum of its full subcategories,
 $${\cA}=
\bigoplus_{D\in \Delta(U,{\Gamma})}{\cA}_{D},$$ where
${\cA}_{D}$ is the restriction of ${\cA}$ on
$D$.

\item\label{enumerate-properties-of-cA-from-finiteness}
If for $a\in U$ holds $\sn\not\in X_{a}(\sm),$ then the class of  $a$ in $\cA(\sm,\sn)$ equals $0.$

\item\label{enumerate-properties-of-cA-to-finiteness}
If for $a\in U$ holds $\sm\not\in X^{a}(\sn),$ then the class of  $a$ in $\cA(\sm,\sn)$ equals $0.$
%
%\item\label{enumerate-properties-of-cA-}
%\item\label{enumerate-properties-of-cA-}
%\item\label{enumerate-properties-of-cA-}
\end{enumerate}
\end{lemma}

\begin{proof}
Let $\sm\in\Specm \Ga$ and $m\geq 1.$
Then $$\lr{\sn^{n}}{\sm^{m}}{U}=\Ga/\sn^{n}\otimes_{\Ga} \r{\sm^{m}}{U}=\Ga/\sn^{n}\otimes_{\Ga} \bigoplus_{l\in X_{m}} \Lr{\sl}{\sm^{m}}{U}=\bigoplus_{l\in X_{m}}\Ga/\sn^{n}\otimes_{\Ga} \Lr{\sl}{\sm^{m}}{U}.$$
If $\sm$ and $\sn$ belong to different classes of $\Delta$-equivalence, then $\sn\ne\sl$ and all the summands in the last sum equal $0.$ It proves \eqref{enumerate-properties-of-cA-disjont-points}. The statement \eqref{enumerate-properties-of-cA-disjoint-cA-D} follows immediately from \eqref{enumerate-properties-of-cA-disjont-points}.
The statements \eqref{enumerate-properties-of-cA-from-finiteness} and \eqref{enumerate-properties-of-cA-to-finiteness} are proved analogously to \eqref{enumerate-properties-of-cA-disjont-points}.
\end{proof}

Let $V$ and $W$ be finitely generated $\Ga$-subbi\-mo\-du\-les in $U$. Then
the bimodule  $T\subset U$, spanned by all products $vw$, $v\in
V,$ $w\in W$, is finitely generated subbimodule in $U$, since $\Ga$
is a Harish-Chandra subalgebra in $U.$  Denote by
$\mu:V\otimes_{\Ga}W\myto T$  the map $\mu(v\otimes
w)=vw,$ $v\in V,w\in W.$

For a  $\Ga$-bimodule $B$  denote by $\lr {\sn^{n}}{\sm^{m}} {\1}
\,  (=\lr {\sn^{n}}{\sm^{m}} {\1}(B))$  the canonical epimorphism
 $$ \lr {\sn^{n}}{\sm^{m}} {\1} :B\myto
 B/(\sn^{n}B+B\sm^{m}).$$

 \begin{lemma}
 \label{lemma-on-the-products-in-cA-approximations}
 Let $V,W\subset U$ be finitely generated $\Ga$-bimodules, $T=VW,$ $\sm,\sn\in \Specm \Ga,$  $m,n\geq 0.$
 \begin{enumerate}
\item\label{enumerate-different-p-th-multiply-as-zero}
 For $\sp, \sp'\in \Specm \Ga,$ $\sp\ne\sp'$ holds $\lR{\sn^{n}}{\sp}{V}\otimes_{\Ga}\Lr{\sp'}{\sm^{n}}{W}=0.$

 \item\label{enumerate-Phi-in-def-cA-is-an-iso}
 The induced by the decomposition \eqref{align-right-left-root-HC-subalgebra-decomposition} homomorphism
 $$\Phi: \bosu{\sp\in X^{\sn}\cap X_{\sm}}{}(\lR{\sn^{n}}{\sp}{V}\otimes_{\Ga}\Lr{\sp}{\sm^{n}}{W})
 \myto \l{\sn^{n}}{V}\otimes_{\Ga}\r{\sm^{m}}{W}$$ is an isomorphism of $\Ga$-bimodules.

\item
\label{enumerate-on-large-enough-middle-term}
There exists $P>0$ ($=P(\sm,\sn,m,n, V,W)$), such that for any $p\geq P$
the canonical projections
$$\lr{\sn^{n}}{\sp^{p}}{\pi} :  \l{\sn^{n}}{V}\myto \lr{\sn^{n}}{\sp^{p}}{V},\, \lr{\sp^{p}}{\sm^{m}}{\pi}:  \r{\sm^{m}}{W}\myto \lr{\sp^{p}}{\sm^{m}}{W}$$
induce isomorphisms
$$ \lR{\sn^{n}}{\sp^{p}} {\pi}:  \lR{\sn^{n}}{\sp}{V}\myto \lr{\sn^{n}}{\sp^{p}}{V},\
\Lr{\sp^{p}}{\sm^{m}}{\pi}:\Lr{\sp}{\sm^{m}}{W}\myto \lr{\sp^{p}}{\sm^{m}}{W}.$$

 \item \label{enumerate-diagram-defined-restricted-from-three-side-composition}

 There exists $P>0,$ such that for $p\geq P$
 there exists a unique homomorphism
 $$ {\widehat{\mu}(\sm^{m},\sn^{n})}:
\bosu{\sp\in X^{\sn}\cap X_{\sm}}{}(\lr{\sn^{n}}{\sp^{p}}{V}
{\otimes}_{\widehat{\Ga}_{\sp}}\lr{\sp^{p}}{\sm^{n}}{W})
\myto \lr{\sn^{n}}{\sm^{m}}{T},$$
which makes the  diagram

\medskip

 \xymatrix @!=.15pc {
V\otimes_{\Ga}W  \ar[rrrrr]^{{}_{\sn^{n}}\1_{\sm^{m}}}\ar[dd]_{\mu}&&&
&&\l{\sn^{n}}{V}\otimes_{\Ga}\r{\sm^{m}}{W}
\ar[dd]_{{}_{\sn^{n}}\mu_{\sm^{m}}}&&&&& & \ar[lllll]_(.7){\Phi}\ar@{->>}[dd]_{\oplus_{\sp} c_{\sp}}\bosu{\sp\in X^{\sn}\cap X_{\sm}}{}(\lR{\sn^{n}}{\sp}{V}
 \otimes_{\Ga}\Lr{\sp}{\sm^{n}}{W})
\\
\\
T  \ar[rrrrr]^{{}_{\sn^{n}}\1_{\sm^{m}} }&&&&&
\lr{\sn^{n}}{\sm^{m}}{T}
&&&&&&
\ar@{=>}[llllll]_(.6){\widehat{\mu}(\sm^{m},\sn^{n})}
\bosu{\sp\in X^{\sn}\cap X_{\sm}}{}(\lr{\sn^{n}}{\sp^{p}}{V}
{\otimes}_{\widehat{\Ga}_{\sp}}\lr{\sp^{p}}{\sm^{n}}{W})
}
\noindent commutative.
Here   ${}_{\sn^{n}}\mu_{\sm^{m}}$ is induced by $\mu,$   $\Phi$ is   defined following \eqref{align-right-left-root-HC-subalgebra-decomposition}
and
 $$c_{\sp}:\lR{\sn^{n}}{\sp}{V}
 \otimes_{\Ga}\Lr{\sp}{\sm^{n}}{W}\xrightarrow{   \quad  \lR{\sn^{n}}{\sp^{p}} {\pi}  \otimes \Lr{\sp^{p}}{\sm^{m}}{\pi}\quad }
 \lr{\sn^{n}}{\sp^{p}}{V}
{\otimes}_{\Ga}\lr{\sp^{p}}{\sm^{n}}{W}
\myto
\lr{\sn^{n}}{\sp^{p}}{V}
{\otimes}_{\widehat{\Ga}_{\sp}}\lr{\sp^{p}}{\sm^{n}}{W},
 $$ where the second arrow is induced by the canonical homomorphism $\Ga\to \widehat{\Ga}_{\sp}.$
 Note that the homomorphism $\widehat{\mu}(\sm^{m},\sn^{n})$ does not depend on $p\geq P.$

 \item\label{enumerate-definition-restricted-composition-by-idempotents}
Denote $$\widehat{\mu}_{\sp}(\sm^{m},\sn^{n}):\lr{\sn^{n}}{\sp^{p}}{V}
{\otimes}_{\widehat{\Ga}_{\sp}}\lr{\sp^{p}}{\sm^{n}}{W}\myto  \lr{\sn^{n}}{\sm^{m}}{T}$$ the restriction  of $\widehat{\mu}(\sm^{m},\sn^{n}).$ Then
$$ \widehat{\mu}(\sm^{m},\sn^{n})=\su{\sp\in X^{\sn}\cap X_{\sm}}{}\widehat{\mu}_{\sp}(\sm^{m},\sn^{n}).$$
Besides, in the assumption of \eqref{enumerate-diagram-defined-restricted-from-three-side-composition}
there exists $S\geq 1,$  such that for every  $\sp\in   X^{\sn}\cap X_{\sm}, $ $v\in V,$ $w\in W$ and its classes $\bar{v}\in \lr{\sn^{n}}{\sp^{p}}{V},$ $\bar{w}\in \lr{\sp^{p}}{\sm^{n}}{W},$ $s\geq S$ and any
\begin{equation}
\label{equation-condition-on-cutting-center}
 \ga\in\Ga, \textrm{ such that  }\ga\equiv 1\mod \sp^{s},   \ga\equiv 0\mod \sq^{s},    \sq\in   X^{\sn}\cap X_{\sm},    \sq\ne \sp
\end{equation}
 holds
$$   \widehat{\mu}_{\sp}(\sm^{m},\sn^{n})(\bar{v}\otimes \bar{w})=\overline{v\ga w},$$
where $\overline{v\ga w}$ is the class of $v\ga w$ in $\lr{\sn^{n}}{\sm^{m}}{T}.$

\item\label{enumerate-correctness-inverse-def-cA-completed-multiplication}
Let $m'\leq m,$ $n'\leq n$ be integers. Then for a sufficiently large  $p$ satisfying the conditions of \eqref{enumerate-on-large-enough-middle-term}, we have the commutative diagram of $\hat{\Ga}_{\sn}-\hat{\Ga}_{\sm}$ homomorphisms
$$
\xymatrix{\lr{\sn^{n}}{\sp^{p}}{V}\otimes_{{\hat{\Ga}_{\sp}}}\lr{\sp^{p}}{\sm^{n}}{W}\ar[d]
 \ar[rr]&& \lr{\sn^{n'}}{\sp^{p}}{V}\otimes_{{\hat{\Ga}_{\sp}}}\lr{\sp^{p}}{\sm^{m'}}{W}\ar[d]
 \\
 \lr{\sn^{n}}{\sm^{m}}{T}
 \ar[rr]&& \lr{\sn^{n'}}{\sm^{m'}}{T},
 }
$$
where the horizontal arrows are induced by the canonical projections and the vertical arrows are    homomorphisms
$\widehat{\mu}_{\sp}(\sm^{m},\sn^{n})$ and $\widehat{\mu}_{\sp}(\sm^{m'},\sn^{n'})$ respectively

\item\label{enumerate-correctness-direct-def-cA-completed-multiplication}
Let $V\subset V',$ $W\subset W'$  be  finitely generated $\Ga$-subbimodules in $U,$ $T'=V'W'.$ Then for a sufficiently large $p$ satisfying the conditions of \eqref{enumerate-on-large-enough-middle-term},  we have the commutative diagram of $\hat{\Ga}_{\sn}- \hat{\Ga}_{\sm}$ homomorphisms
$$
\xymatrix{\lr{\sn^{n}}{\sp^{p}}{V}\otimes_{{\hat{\Ga}_{\sp}}}\lr{\sp^{p}}{\sm^{n}}{W}\ar[d]
 \ar[rr]&& \lr{\sn^{n}}{\sp^{p}}{V'}\otimes_{{\hat{\Ga}_{\sp}}}\lr{\sp^{p}}{\sm^{m}}{W'}\ar[d]
 \\
 \lr{\sn^{n}}{\sp^{p}}{T}
 \ar[rr]&& \lr{\sn^{n}}{\sp^{p}}{T'},
 }
$$
where the horizontal arrows are induced by the canonical embeddings and the vertical arrows are  corresponding homomorphisms
$\widehat{\mu}_{\sp}(\sm^{m},\sn^{n}).$

 \end{enumerate}
 \end{lemma}

 \begin{proof}
 The statement \eqref{enumerate-different-p-th-multiply-as-zero} obviously follows from the Chinese remainder theorem. The statement \eqref{enumerate-Phi-in-def-cA-is-an-iso} follows from the decomposition \eqref{align-right-left-root-HC-subalgebra-decomposition} and from the statement \eqref{enumerate-different-p-th-multiply-as-zero} above.
To prove the statement \eqref{enumerate-on-large-enough-middle-term} note first that the sequence of finite dimensional modules  $\lr{\sn^{n}}{\sp^{p}}{V}$ stabilizes for some $P$ and all $p\geq P$. Hence   $\lR{\sn^{n}}{\sp}{V}$ is a factor of $\lr{\sn^{n}}{\sp^{p}}{V}$, $p\geq P$. On the other hand every $\lr{\sn^{n}}{\sp^{p}}{\pi}$ factorizes through  $\lR{\sn^{n}}{\sp}{V},$ which proves \eqref{enumerate-on-large-enough-middle-term} for $V$. The case of $W$ is considered analogously.

 Note  that for sufficiently large $p$ all $c_{\sp}$ are isomorphisms. In fact, the first map in the definition of $c_{\sp}$ is an isomorphism due to \eqref{enumerate-on-large-enough-middle-term}, and the second map is an isomorphism, since both multipliers in the tensor product are finite dimensional $\sp$-torsion modules. Hence the third  vertical arrow  in the diagram is an isomorphism. Besides, $\Phi$ is an isomorphism by \eqref{enumerate-Phi-in-def-cA-is-an-iso}, that proves \eqref{enumerate-diagram-defined-restricted-from-three-side-composition}. Independence on $p$ is obvious.

 Let $\bar{v}=\su{\sl}{}{\r{\sl}{\bar{v}}},$ $\bar{w}=\su{\sl}{}{\l{\sl}{\bar{w}}},$   be the decompositions from \eqref{align-right-left-root-HC-subalgebra-decomposition}. Then for $\ga',\ga''\in\Ga,$ satisfying \eqref{equation-condition-on-cutting-center} for large enough $s$ holds
 $$ \bar{v}\ga'=\r {\sp}{\bar{v}}+\su{\sl\not\in  X^{\sn}\cap X_{\sm}}{}\r{\sl}{\bar{v}}\ga',\quad \ga''\bar{w}=\l{\sp}{\bar{w}}+\su{\sl\not\in  X^{\sn}\cap X_{\sm}}{}\ga''\l{\sl}{\bar{w}}.$$
 Then following \eqref{enumerate-different-p-th-multiply-as-zero}, $v\ga'\otimes \ga'' w=\r{\sp}{v}\otimes \l{\sp}{w},$ which competes the proof of \eqref{enumerate-definition-restricted-composition-by-idempotents}.
 The statements \eqref{enumerate-correctness-inverse-def-cA-completed-multiplication} and \eqref{enumerate-correctness-direct-def-cA-completed-multiplication} follow from \eqref{enumerate-definition-restricted-composition-by-idempotents}. In these diagrams  both images of $\bar{v}\otimes\bar{w}$ from the left upper corner of the diagram   equal to the class $\overline{v\ga w}$ in the right lower corner for a suitable $\ga\in\Ga$. \end{proof}

A composition in the category $\cA$ is defined as follows.
Since direct limits commute with tensor product, we can write

\begin{equation}
\label{equation-def-of-composition}
\cA(\sl,\sn)\otimes_{\widehat{\Ga}_{\sl}} \cA(\sm,\sl)\simeq \varinjlim_{V\in F(U)}\varinjlim_{W\in F(U)} \lr{\sn}{\sl}{V}
\otimes_{\widehat{\Ga}_{\sl}} \lr{\sl}{\sm}{W}
\end{equation}

Then we have the following composition
\begin{equation*}
\lr{\sn}{\sl}{V}
\otimes_{\widehat{\Ga}_{\sl}} \lr{\sl}{\sm}{W}\xrightarrow{ \ {\widehat{\mu}_{\sl}(\sm^{m},\sn^{n})} \ } \lr{\sn^{n}}{\sm^{m}}{T} \xrightarrow{\  i_{\sn,\sm} \ } \cA(\sm,\sn),
\end{equation*}
where the first homomorphism is constructed above and the second is the canonical map in the direct limit.  From the commutative diagrams in
 Lemma \ref{lemma-on-the-products-in-cA-approximations}, \eqref{enumerate-correctness-inverse-def-cA-completed-multiplication}, we have a well defined map

 $$\cA(\sp,\sn)\times \cA(\sm,\sp)\myto \cA(\sm,\sn).$$

\subsection{Generalized Pierce decomposition}
\label{Generalized-Pierce-decomposition}
 We start from the
following categorical statement. Assume $\cC$ is  a $\k$-category
with sums and products and $\{C_{i}\mid i\in \cI\}$ be a family of
objects of $\cC$. Denote by (*)   the following properties
of this family:

(*) for every $j\in \cI$ there exist finitely  many $i\in \cI$, such that $\cC(C_{i},C_{j})$ and $\cC(C_{j},C_{i})$ are nonzero.

Consider the vector space $$\Pi_{\cI}=\prod_{(i ,j)\in\cI\times \cI} \cC(C_{j},C_{i}),$$ written as $\cI\times \cI$ matrices, provided that $j$'s correspond to columns and $i$'s corresponds to rows. In general the standard "column-by-row" product of such matrices is not well defined.
By $\rM_{\cI}$ denote the subspaces of $\Pi_{\cI},$ formed by the matrices with finitely many nonzero elements in any column  and in any row. Then   "column-by-row" product turns it into a $\k$-algebra.

\begin{lemma}
\label{lemma-cA-as-endomorphisms-of-direct-sum}
Assume the family $\{C_{i}\mid i\in \cI\}$ of objects of the $k$-category $\cC$ satisfies the property (*).
\begin{enumerate}
\item\label{enumerate-from-category-to-matrix}
There exists a canonical isomorphism of $\k$-algebras
$$\rM_{\cI} \simeq \End_{\cC}(\osu{i\in \cI}{} C_{i}),$$
where  $\End_{\cC}$ denotes the endomorphisms ring in the category $\cC.$
\item\label{enumerate-from-category-to-matrix-equivalence}
Let $\cC_{\cI}$ be the restriction of the category $\cC$ on $\cI$,
 $\rM_{\cI}-\Mod_{r}$  the full subcategory in $\rM_{\cI}-\Mod$ consisting of modules $M,$ such that $M=\bosu{i\in \cI}{} e_{ii} M,$  where $e_{ii}$ is a matrix unit corresponding  $i\in \cI.$
Then there exists the canonical equivalence
$$F: \cC_{\cI}-\Mod\simeq \rM_{\cI}-\Mod_{r},$$   where for
$N\in  \cC_{\cI}-\Mod$  holds $F(N)=\bosu{i\in\cI}{}N(i).
$

\end{enumerate}
\end{lemma}

\begin{proof}
Every element $(f_{ij}\mid i,j\in \cI)\in \Pi_{\cI}$ defines canonically a homomorphism $f:\osu{  i\in \cI}{}  C_{i}\myto \pru{ i\in \cI} {}C_{i}.$ By the condition (*) the image of $f$ belongs to $\osu{  i\in \cI}{}  C_{i}\subset \pru{ i\in \cI} {}C_{i}.$
The second statement is standard.
\end{proof}

The following statement is an analogue of two-sided Pierce decomposition for the pair $\Ga\subset U.$

\begin{theorem}
\label{theorem-embedding-in-matrix-ring}
For $u\in U$ denote  by $[u]$ the matrix from $ \rM_{\cA},$ such that  $[u]_{\sm,\sn} =(u_{\sn,\sm}), $ $\sm,\sn\in\Specm \Ga,$ where $u_{\sn,\sm}$ is the image of $u$ in $\cA(\sm,\sn).$
\begin{enumerate}
\item\label{enumerate-embedding-is-homomorphism}
The mapping $i_{\cA}=[\,]:U\myto \rM_{\cA}$, which sends $u\in U$ to  $[u]$ is a homomorphism of $\k$-algebras.
\item\label{enumerate-embedding--image-is-dense}
Endow $\Pi_{\cA}$ with the topology of direct product and $\rM_{D}\subset \Pi_{\cA}$ with the induced topology.
Then the image of $\pi_{D}i_{\cA}$ is dense in $\rM_{D}.$
\end{enumerate}
\end{theorem}

\begin{proof}
Following Lemma \ref{lemma-simplest-properties-of-cA},
\eqref{enumerate-properties-of-cA-from-finiteness},
\eqref{enumerate-properties-of-cA-to-finiteness}, the matrix $[u]$
has finitely many nonzero elements in any row and any column,
hence $[u]\in\rM_{\cA}.$

Fix $u,v\in U,$ $\sm,\sn\in \Specm \Ga,$ $m,n\geq1$. Also fix
$V,W,T$ and  $P>0,$ satisfying the conditions of Lemma
\ref{lemma-on-the-products-in-cA-approximations}. Then from
 $$ \widehat{\mu}(\sm^{m},\sn^{n})=\su{\sp\in X^{\sn}\cap X_{\sm}}{}\widehat{\mu}_{\sp}(\sm^{m},\sn^{n})$$ follows
 (as in the proof of Lemma \ref{lemma-on-the-products-in-cA-approximations},  \eqref{enumerate-definition-restricted-composition-by-idempotents})
 $$\overline{uv}= \su{\sp\in  X^{\sn}\cap X_{\sm}}{}\widehat{\mu}_{\sp}(\sm^{m},\sn^{n})(\r{\sp}{\bar{v}}\otimes  \l{\sp}{\bar{w}}),$$
where $\overline{uv}$ is a class of $uv \in \lr{\sn^{n}}{\sm^{m}}{T}.$ Taking the limits we obtain that
$$ [uv]_{\sm,\sn}=  \su{\sp\in  X^{\sn}\cap X_{\sm}}{} [u]_{\sn,\sp} [v]_{\sp,\sm}, $$
that proves the first statement.
To prove the second statement it is enough to show that if
$\cA(\sm,\sn)\subset \rM_{D},$ then  $\cA(\sm,\sn)\subset
\overline{\Im\pi_{D}i_{\cA}}.$ First note that the image of $U$ is
dense in the image of $\cA(\sm,\sn)\subset\rM_{D},$ that is, matrices from $\rM_{D},$ which are zero in all positions except  $(\sn,\sm).$ Let $\{u_{i}\in U\mid
i=1,\dots \}$ be a sequence in $U$, which converges to an element
$f$ from $\cA(\sm,\sn).$
Note that by definition any $D$ is  at
most countable. Consider an increasing sequence of
finite subsets $S_{i}\subset D,$ $i=1,\dots$  such that $\cu{i=1}{\infty} S_{i}=D,$  strictly increasing sequence of integers $k_{i}>0$ and the elements
$\mu_{i},\nu_{i}\in \Ga,$ such that

\begin{equation*}
\left\{\begin{array}{cccccc}
  \mu_{i} & \equiv & 1 \ \mod & \sm^{k_{i}} \\
  \mu_{i} & \equiv & 0 \ \mod & {\sm'}^{k_{i}}, &\sm'\in S_{i},\sm'\ne\sm.
\end{array}\right.
\end{equation*}
and
\begin{equation*}
\left\{\begin{array}{cccccc}
  \nu_{i} & \equiv & 1  \mod & \sn^{k_{i}} \\
  \nu_{i} & \equiv & 0  \mod & {\sn'}^{k_{i}}, &\sn'\in S_{i},\sn'\ne\sn.
\end{array}\right.
\end{equation*}

Then $i_{D}(\mu_{i})$ (respectively $i_{D}(\nu_{i})$) tends in $\rM_{D}$ to the diagonal matrix unit in the position $\sm$ (respectively $\sn$), hence   the sequence $i_{D}(\nu_{i}u_{i}\mu_{i})=$ $i_{D}(u_{i})$ $i_{D}(\mu_{i})$ converges to $f$ since it tends to $0$ in all positions except $(\sn,\sm)$ and tends to $f$ in the position $(\sn,\sm).$

\end{proof}

\subsection{Gelfand-Tsetlin modules as $\cA$-modules}
\label{subsection-Gelfand-Tsetlin-modules-asScAS-modules}
We consider  the category of ${\k}-{\Mod}$  endowed with the
discrete topology and consider the category ${\cA}-\Mod_d$ of continuous functors $M:{\cA}{\longrightarrow}{\k}-\Mod,$ which in \cite{dfo:hc} are called \emph{discrete modules}{}. It means, for every $\sm\in\Specm \Ga$ there exists $m (=m(\sm)) \geq0,$ such that $\sm^{m}M(\sm)=0.$
For a discrete ${\cA}-$module $M$ define a
Gelfand-Tsetlin $U-$module $\bF(M)$ by setting
\begin{align}
\label{align-definition-of-action-on-module}
&\bF(M)=\bosu{\sm\in\Specm {\Gamma}}{} M(\sm) \text{ and for }  x\in
M(\sm), a\in U\\ \nonumber & \text{ set } ax= \sum_{\sn\in \Specm {\Gamma}} a_{\sm,\sn}x,
\end{align}
where
$a_{\sm,\sn}$ is the image of $a$ in ${\cA}(\sm,\sn)$. If $f:M{\longrightarrow} N$
is a morphism in ${\cA}-{\mo}_d$ then define
$\bF(f)=$ $\bosu{\sm\in \Specm {\Gamma}}{}
f({\sm})$.
\begin{theorem}\label{theorem-HC-mod-equiv-A-mod} (\cite{dfo:hc}, Theorem 17)
The defined above functor $\bF$ is an equivalence of categories
$$\bF:{\cA}-{\mo}_d {\longrightarrow}
\bH(U,{\Gamma}).$$ Moreover it
  induces a functorial isomorphism
$$\Ext_{\cA}^{1}(\mathbb{F}(X),\mathbb{F}(Y))\simeq\Ext_{U}^{1}(X,Y), \ X,Y\in \bH W(U,{\Gamma}).$$
\end{theorem}

\begin{proof}
Following Lemma \ref{lemma-cA-as-endomorphisms-of-direct-sum}, \eqref{enumerate-from-category-to-matrix-equivalence} there exists a canonical equivalence
$$ \cA-\mo_{d}\simeq \rM_{\cA}-\mo_{r},$$  that together with the functor, induced by  $i_{\cA}:U\hookrightarrow  \rM_{\cA}$ (Theorem \ref{theorem-embedding-in-matrix-ring}, \eqref{enumerate-embedding-is-homomorphism}) gives us the functor $\cA-\mo_{d}\myto U-\mo,$ besides the image of this functor belongs to $\bH(U,{\Gamma}).$
The statement on $\Ext^{1}$ follows from the fact that the functor $\bF$ preserves the exact sequences.
\end{proof}

Let $u\in U$, $\sm,\sn\in \Specm \Ga$,
${u}_{\sm,\sn}$ the image of $u$ in
$\cA(\sm,\sn)$. Let $$M=\underset{\sm\in\Specm \Ga}\sum M(\sm)$$
be a Gelfand-Tsetlin module. The action of
${u}_{\sm,\sn}$ on $M$ is described in the next lemma.

\begin{lemma}
\label{lemma-action-of-the-class} Let $m_{u}:M(\sm)\myto M$
be the map $x\longmapsto ux$, $x\in M(\sm)$,
$$u_{\sm,\sn}=\pi_{\sn}m_{u}: M(\sm)\to M(\sn),$$
where $\pi_{\sn}$ is the canonical projection of $M$ onto
$M_{\sn}$. If $\sn\not\in X_{u}(\sm)$, then
$u_{\sm,\sn}=0$.
\end{lemma}
\begin{proof} Condition $\sn\not\in X_{u}(\sm)$ holds if and only if $\Ga u
 \Ga =\sn u\Ga+\Ga u\sm$, equivalently, $u\in \sn u\Ga+\Ga
 u\sm$. But then $u\in \sn^{n} u\Ga+\Ga
 u\sm^{m} $ for every $m,n>0$ and hence $u_{\sm,\sn}=0$.
On the other hand in this case $\pi_{\sn}m_{u}=0$ by
\eqref{equation-1}.
By the Chinese remainder theorem, there exists a sequence
$\gamma_{n}\in\Gamma,n\geq 1$, such that $\gamma_{n}\equiv
1 (\mod \sn^{n})$, $\gamma_{n}\equiv0 (\mod \sl^{n})$ for any
$\sl \in X_{u}(\sm),\sl\ne\sn$. Then in $\cA(\sm,\sn)$
holds $\ds\lim_{n\to\infty}\gamma_{n}u={u}_{\sm,\sn}$.
On the other hand the operator $\bar{m}_{u}=$
$\ds\lim_{n\to\infty} m_{\gamma_{n}u}:M(\sm)\to M$ is
well defined,  $\pi_{\sn}m_{u}=\pi_{\sm}\bar{m}_{u}$
and $\pi_{\sn'}m_{u}=0$ for any $\sn'\ne \sn$. Then
\eqref{align-definition-of-action-on-module} completes
the proof.
\end{proof}

We will identify a discrete ${\cA}$-module $N$ with
the corresponding Gelfand-Tsetlin module $\bF(N)$. For $\sm\in
\Specm \Gamma$ denote by $\hat\sm$ a completion of
$\sm$. Consider a two-sided ideal $I\subset \cA$ generated by $\hat\sm$ for all $\sm\in
\Specm \Gamma$ and set ${\cA}(W)= {\cA}/I$.
Then Theorem \ref{theorem-HC-mod-equiv-A-mod} implies the
following corollary.

\begin{corollary}
\label{cor:hc-weight-equiv-Aw} The categories $ \bH W(U,{\Gamma})$  and ${\cA}(W)-{\mo}_{d}$ are
equivalent.
\end{corollary}

The following standard statement shows usefulness of the category $\cA$ for the study of simple Gelfand-Tsetlin modules over $U$.

\begin{lemma}
\label{lemma-simple-modules-in-vertex}
Let  $M$ be a simple $\cA$-module, $\sm\in \Specm \Gamma$,
 $M(\sm)\ne0$, $N$ a simple $\cA(\sm,\sm)$-module. Then the correspondences
\begin{align}
\label{align-simple-modules-in-vertex}\nonumber
&M\longmapsto M|_{\cA(\sm,\sm)} \text{ and } \
N\longmapsto (\cA\otimes_{\cA(\sm,\sm)}N)/J,
\end{align}
 where  $J$  is a unique maximal $\cA$-submodule in the induced module,
realizes a bijection between the sets of isomorphism classes of simple $\cA(\sm,\sm)$-modules  and isomorphism classes of simple $\cA$-modules $M$ such that $M(\sm)\ne0$.
\end{lemma}

The subalgebra ${\Gamma}$ is called  \emph {big}{} in
$\sm\in$ $\Specm {\Gamma}$ if ${\cA}(\sm,\sm)$ is finitely generated as a right
${\Gamma}_{\sm}-$module. The importance of this
concept is described in the following statement.

\begin{lemma}\label{lemma-Gamma-big} (\cite{dfo:hc}, Corollary 19) If
${\Gamma}$ is big in $\sm\in$ $\Specm {\Gamma}$ then there
exists finitely many non-isomorphic irreducible Gelfand-Tsetlin
$U-$modules $M$ such that $M(\sm)\ne 0$. For any such module
$M$ holds $\dim_{\k} M(\sm)<\infty$.
\end{lemma}

\subsection{Examples of computation of $\cA$}
\label{subsection-Some-examples-of-calculation-of-cA} First
example is given in the beginning of the section, namely, the
presentation of a basic associative algebra as a quiver with
relations.

Now we illustrate our techniques applying it to the study of
representations of skew group algebras over a commutative ring and
to obtain the classical result on the construction of irreducible
representations of a finite group $G=N\rtimes H$ with abelian $N.$

The case of a skew group algebra is summmarized  in the
following statement.

\begin{proposition}
\label{proposition-HC-modules-for-skew-group-algebras}
Let $\La$ be a commutative ring, $\cM$  a monoid, acting on $\La,$
$U=\La*\cM$ the skew-group algebras.  For $\sm,\sn \in\Specm \La$ set
$$\cM(\sm,\sn)=\{\vi\in\cM\mid \vi\cdot \sm=\sn\}.$$ Then $\La\subset U$ is a Harish-Chandra subalgebra and

\begin{enumerate}
\item\label{enumerate-hc-skew-group-objects-of-blocs} The blocks
of the category $\cA=\cA(\La)$ correspond to the orbits $\Spec
\La/\cM.$ For $D\in \Spec \La/\cM$ the set of objects of $\cA_{D}$
coincides with $D$. Moreover, $\cA_{D}$ itself is a groupoid  such
that for $\sm\in D$ holds $\cA(\sm,\sm)\simeq
\La_{\sm}*\cM(\sm,\sm).$

\item\label{enumerate-hc-skew-group-existence-of-irreducible-modules}
In every block there exists a simple Gelfand-Tsetlin module.

\item\label{enumerate-hc-skew-group-free-action}
If the action of $\cM$  on $D$ is free, then $\cA_{D}-\mo_{d}$ is equivalent to the category of finite dimensional
modules over $\La_{\sm}$ for any $\sm\in \Specm \La.$
%\item\label{enumerate-hc-skew-group-}
\end{enumerate}
\end{proposition}

\begin{proof}
To compute $\cA(\sm,\sn)$ there is enough to consider in
Definition \ref{definition-of-category-cA}  the bimodules $V$ of
the form $V=\bosu{\vi\in S}{}\La\vi$ for some $S\subset \cM.$ In
this case
\begin{equation*}
V/(\sn^{n}V+V\sm^{m})\simeq \bosu{\vi\in S\cap\cM(\sm,\sn) }{}(\La/\sn^{k})\vi, \text{ where }k=\min\{m,n\}.
\end{equation*}
Since in the definition of $\cA(\sm,\sn)$ we can consider only $V$ containing $\cM(\sm,\sn),$ we obtain that
\begin{equation}
\label{equation-cA-in-skew-group}
\cA(\sm,\sn)\simeq \bosu{\vi\in \cM(\sm,\sn) }{}\La_{\sn}\vi \simeq \bosu{\vi\in \cM(\sm,\sn) }{}\vi\La_{\sm}.
\end{equation}
The last isomorphism ($\La_{\sn}\vi \simeq \vi \La_{\sm}$) is
defined using the isomorphism of $\La$-bimodules
$$(\La/\sn^{k})\vi\simeq \vi(\La/\sm^{k}), k\geq 1. $$ This
isomorphism allows to define the composition in $\cA$ as in
Theorem \ref{theorem-embedding-in-matrix-ring}. In particular
$\cA(\sm,\sm)$  is isomorphic to $\bosu{\vi\in \cM(\sm,\sm)
}{}\La_{\sm}\vi$ as $\La$-bimodule.  It  is also isomorphic to
$\La_{\sm}*\cM(\sm,\sm)$ as an algebra. Any pair $(\sm,\sn)$ of
the objects in a block $\cA$ are isomorphic by application of an
element from $\cM(\sm,\sn),$ which proves
\eqref{enumerate-hc-skew-group-objects-of-blocs}.

By Lemma \ref{lemma-simple-modules-in-vertex} the statement
\eqref{enumerate-hc-skew-group-objects-of-blocs} reduces the
problem of classification of simple Harish-Chandra $U$-modules to
the problem of classification of simple
$\La_{\sm}*\cM(\sm,\sm)$-modules, $\sm\in\Specm \La.$ But the
$\La$-bimodule $\La_{\sm}$ is a direct summand of
$\La_{\sm}*\cM(\sm,\sm)$ as a $\La_{\sm}$-bimodule, hence
$U/U\sm\ne0.$ It proves
\eqref{enumerate-hc-skew-group-existence-of-irreducible-modules}.

The statement \eqref{enumerate-hc-skew-group-free-action} is obvious.
\end{proof}

Next   we will show that the classical Mackey construction can be
interpreted as a special case of the category of Gelfand-Tsetlin
modules (\cite{fo-revista}). We assume that the base field $\k$ is
algebraically closed and its characteristic is coprime with
$n=|G|.$ We set $U=\k[G]$ and $\Gamma=\k[N]$. Obviously, $\Gamma$
is a Harish-Chandra subalgebra in $U$. Denote by $\check{N}$ the
set of characters of $N$. Then $\k[N]\simeq
\ds\prod_{\chi\in\check{N} }\k_{\chi}$, where $\k_{\chi}$
corresponds to $\chi$.  By $\sm_{\chi}$ denote the kernel the
character $\chi$.

The group $G$ acts by conjugations from the left on the group $N$,
$x\longmapsto {}^{g}x=gxg^{-1}, x\in N,g\in G$, and it also acts
on $\check{N}$ from the right, $\chi\longmapsto {\chi}^{g},
\chi\in \check{N},g\in G$, such that $\chi^{g}(x)=
\chi({}^{g}x),g\in G$. Denote by $\St_{G} x$ and $\St_{G}\chi$ the
stabilizers of the corresponding actions.

We give a construction of the category $\cA=\cA_{U,\Gamma}$. In
\cite{fo-revista}    the following statement is shown:

\begin{proposition}
\label{proposition-mackey-gruppoid} Let $\kY=\kY(G,N)$ be a
groupoid, such that $$\Ob \kY=\check{N}, \quad
\kY(\chi_{1},\chi_{2})=\{g\in G\,|\, \chi_{1}=\chi_{2}^{g}\},
\quad \chi_{1},\chi_{2}\in \check{N},$$ $\kN$ be a subgroupoid of
$\kY$, such that $\kN(\chi, \chi)=N$ for any $\chi\in\check{N}$
and empty otherwise, $\kX=\kY/\kN$ and $\k\kX$ be its $\k$-linear
envelope. Then there exists a canonical isomorphism of categories
$\Phi:\k\kX\myto \cA,$ identical on the objects.
Besides
\begin{itemize}
\item[1.] \label{step-endomorphism-ring} $\cA(\chi,\chi)\simeq
\k[\St_{\chi}]$.
\item[2.]\label{step-isomorphic-characters}$\chi_{i},\chi_{j}\in\Ob
\cA$ are isomorphic if and only if $\chi_{i}^{g}=\chi_{j}$ for some $g\in G.$
\end{itemize}
\end{proposition}

\section{Generic representations}
\label{section-Representations-of-general-position} In this
section we will assume that $\cM$ is finitely generated, in
particular, $\cM$ is countable. Assume $\k$ is uncountable. A
non-empty set $X\subset \Specm \Ga$ is called \emph{massive}{},
provided that $X$ contains an intersection of countable many dense
open subsets.  In this case $X$ is dense in $\Specm \Ga$
(\cite{ga}, Lemma 4.2). Assume $s:\Ga\myto \Ga'$ is a homomorphism
of algebras, such that for the induced $s^{*}:\Specm
\Ga'\myto\Specm \Ga$, $\Im s^{*}$ is massive in $\Specm \Ga.$ In
this case if $\cX'\subset \Specm \Ga'$ is a massive subset, then
$s^{*}(\cX')$ is a massive subset in $\Specm \Ga.$

Let $U_{\chi}=U/(U \sm_{\chi})$. We call $U_{\chi}$ the
\emph{universal module generated by a $\chi$-eigen\-vec\-tor}. If $\Ga$ is a Harish-Chandra subalgebra then $U_{\chi}\in $ $\bH(U,{\Gamma})$. The character $\chi$ lifts to a simple $U$-module if and only if $U_{\chi}$ is nonzero.

Let $\La_{0}$ be the algebra generated by $\underset{m\in\cM}\bigcup
m(\Gamma)$ (if ${\Ga}_{0}\subset \bar{\Gamma}$ then $U$ is a right Galois order). Every nonzero element $u\in U$ can be presented (not uniquely) in the form
$$u=\sum_{\vi\in\cM/G}[a(u)_{\vi}\vi], \textrm{ where } a(u)_{\vi}\ne 0. $$
Denote by $\La_{1}$ the algebra, generated by $\La_{0}$ and all possible $m\cdot a(u)_{\vi},$ where $u$ runs $U\setminus\{0\}$ and and $m$ runs $\cM.$ Analogously by $\La_{2}$ the algebra generated by $\La_{0}$ and all $(m a(u)_{\vi})^{\pm 1}.$ Also denote $\cL_{i}=\Specm \La_{i}$, $\pi_{0}$ the canonical map from $\cL_{0}$ to $\Specm \Ga,$ $\Om_{i}$ $=$ $\pi_0(\cL_{i})$, $i=1,2$. Then
$$\La_{2}\supset\La_{1}\supset\La_{0},\ \cL_{2}\subset \cL_{1}\subset\cL_{}, \,\,\,\,\,
 \Om_{2}\subset \Om_{1}\subset\Om.$$
By $\cL_{r}\subset\cL$ denote the set of $\ell$, such that $\cM$ acts on $\ell$ freely and $\cM\cdot\ell\cap G\cdot\ell=\{\ell\}$ (in terms of \ref{subsection-Finiteness-of-extensions-for-integral-Galois-algebras}
for $\sm=\pi_0(\ell)\in \Specm \Ga$ holds $S(\sm,\sm)=\{e\}$) and
$\Om_{r}=\pi_0(\cL_{r})$.

%
%If $\cM\cdot \Ga\subset \bar{\Ga}$ then any $m\in\cM$ defines an
%automorphism of $\bar{\Ga}_0=\bar{\Ga}$ and, hence induces an
%automorphism of $\cL$. In particular, this holds when $\Ga$ is a
%Harish-Chandra subalgebra.

\begin{lemma}\label{lemma-regular-tableau-massive}
The sets $\cL_{i}\subset\cL$ and $\Om_{i}\subset \Specm \Ga$,
$i=1,2,$ are massive. Besides $\cL_{i},$ $i=1,2$ are invariant with respect to the actions of $G$ and $\cM.$ If $\cM\cdot \Ga\subset
\bar{\Ga}$ then $\cL_r$ and $\Omega_r$ are massive.
\end{lemma}
\begin{proof}
Every element $a\in L$ defines a rational function on $\cL$ with a
non-empty open domain of  $X(a)\subset \cL$. Then
$\cL_{i},i=1,2$ by definition are intersections of countably many
sets of the form $X(a)$, hence are massive. The properties of $\cM$- and $G$-invariance follows from the facts, that $gm\Ga=m^{g}\Ga$ and
$[a\vi]=[a^{g}\vi^{g}],$ $g\in G,$ $m\in\cM.$

For any $m\in \cM, m\ne e$, denote $\cX_{m}= \{\ell\in\cL \mid m
\cdot \ell\in G\cdot \ell\}$ and for $m\in \cM, g\in G$ define a
closed in $\cL$ subset $\cL(m,g)=\{\ell\in \cL\mid m\cdot
\ell=g\cdot \ell\}$. Then $\cX_{m}=\underset{g\in G}\cup\cL(m,g)$
is closed subset in $\cL$, since $G$ is finite. Assume that
$\cL=\cX_{m}$ for some $m\ne e$.  Since ${\cL}$ is irreducible, we
conclude that $\cL(m,g)=\cL$ for some $g\in G$, and hence $m=g$.
This is a contradiction  since $\cM$ is separating. But $\underset{m\in\cM,m\ne
e}\cup\cX_{m}$ is the complement of ${\\cL_{r}}$ in $\cL$ and
$\cL_r$ is massive. The sets $\Omega_i$, $i=1,2,r$, are massive as the images of massive subset $\cL_{i}.$
\end{proof}

The following useful fact follows from the separating of $\cM$-action (Lemma \ref{lemma-separ-monoid}).
\begin{lemma}
\label{lemma-distinguish-orbits} Let $\bar{\Gamma}$ be the
integral closure on $\Ga$ in $L$, $\ell_{\sm_{k}},\dots,\ell_{\sm_{k}}\in\Specm
\bar{\Gamma}$ belong to different orbits of $G$. Then there
exists $\gamma\in \Gamma$, such that
$\gamma(\ell_{\sm_{1}}),\dots,\gamma(\ell_{\sm_{k}})$ are distinct, that is
$\Gamma$ distinguishes the orbits of $G$.
\end{lemma}

\begin{theorem}
\label{proposition-repr-of-schat} Suppose that the field
$\k$ is uncountable, $\chi$ a character of $\Ga,$ $\sm_{\chi}\in\Specm \Ga$ its kernel.

\begin{enumerate}
\item\label{enum-massive-X0} If $\chi\in \Omega_1$, then $U_{\chi}$ is nonzero
    Gelfand-Tsetlin module and $\supp
    U_{\chi}\subset\cM\sm.$

\item\label{enum-massive-X1} If $\cM$ is a group, then for any $\chi\in X_{2}$ the module
    $U_{\chi}$ is a unique $\cO_{\sm}$-graded
    irreducible $U$-module generated by a
    $\chi$-eigenvector and $\supp U_{\chi}=\cO_{\sm}$.

\item\label{enum-nonperiodic} If $\cM$ is a group and
    $\cM\cdot \Ga\subset \bar{\Ga}$,  $\chi\in \Omega_{2}\cap \Omega_{r}$, then  $U_{\chi}$ is irreducible weight
$U$-module with   $1$-dimensional components. In this case there
is a canonical isomorphism of $\k$-vector spaces $\k\cM\simeq
U_\chi$.
\end{enumerate}
\end{theorem}

\begin{proof}
Note that the canonical embedding $U\hookrightarrow L*\cM$ factorizes through a skew semigroup algebra $\Lambda_1*\cM\subset L*\cM.$ A character $\tilde{\chi}:\Lambda_1\rightarrow \k$
induces a $\Lambda_1*\cM$-module $M(\tilde{\chi})$ which is
isomorphic to $\underset{\vi \in \cM}
\oplus\Lambda_{1}/(\tilde{\sm}^{\vi}), $ as a $\k$-vector space,
where $\tilde{\sm}=\Ker \tilde{\chi}$.

Since the
restriction of $M(\tilde{\chi})$ on $U$ is nontrivial then for every
$\chi\in X_{1}$, $U_{\chi}\neq 0$. Clearly, $U_{\chi}$ is a
Gelfand-Tsetlin
 module.
Moreover, since $U_{\chi}$ is $\Specm \cL$-graded, it has an
irreducible quotient with a nonzero $\chi$-eigenvector. This
implies \eqref{enum-massive-X0}.

 Let $\chi\in \Omega_{2}$,  $u=\su{h\in
\cM}{} [x_{h} h]\in U$. By assumption, for every
$\sn\in \cM\cdot \sm $ holds $x_{h}(\sn)\ne 0$, hence every graded
component of $U_{\chi}$  generates the
whole $U_{\chi}$. Therefore, $U_{\chi}$ is irreducible as
$\cO_{\sm}$-graded $U$-module. Moreover, since $U_{\chi}$
is the universal module generated by a $\chi$-eigenvector,
it is a unique such graded irreducible module, implying
statement \eqref{enum-massive-X1}.

Note that if $\cM$ acts
on $\sm$ with a nontrivial stabilizer then $U_{\chi}$ is
not irreducible. Since $\cM\cdot \Ga\subset
\bar{\Ga}$, the set $\Omega_r$ is massive by
Lemma~\ref{lemma-regular-tableau-massive}. Consider a
subset $\Omega_2\cap \Omega_r$. Since $\Gamma$ distinguishes
$G$-orbits by Lemma~\ref{lemma-distinguish-orbits}, it
implies the irreducibility of $U_{\chi}$ for any $\chi\in
X_r$. The basis elements of $U_{\chi}$ in this case are
labeled by the elements of $\k\cM$ which completes the
proof of \eqref{enum-nonperiodic}.
\end{proof}

\begin{lemma}
\label{lemma-tor-is-nonzero-on-closed-set} Let  $M$  be a
finitely generated right $\Ga$-module. Then the set of
$\sm\in\Specm \Ga$ such that $\Tor_{1}^{\Ga}(M,\Ga/\sm)=0$
contains an open dense subset  $X\subset\Specm \Ga$.
\end{lemma}
\begin{proof}
Let $ R^{\bullet}\ : \ \dots\xar{\ d^{2}} \Ga^{n_{2}}\xar{\
d^{1}}\Ga^{n_{1}}\xar{\ d^{0}}\Ga^{n_{0}} \longrightarrow
0\dots$  be a free resolution of $M$. It induces the
resolution $R^{\bullet}\otimes_{\Ga}K$ of
$M\otimes_{\Ga}K$.  Denote $r=\dim_{K} \Im(d^{1}\otimes
\id_{K})$. Let $d^{i}(\sm)$ be the specialization of $d^{i}$
in $\sm\in\Specm \Ga$. Then all those $\sm$ satisfying the
conditions $\dim \Im d^{1}(\sm)= n_{1}-r$ and $\dim \Im
d^{2}(\sm)= r$, form an open dense set $X$. Then
for $\sm\in X$ the first cohomology of
$R^{\bullet}\otimes_{\Ga} \Ga/\sm$ equals $0$, that completes the proof.
\end{proof}

\begin{proposition}
\label{proposition-embedding-U-in-cA} For a nonzero $u\in
U,u\ne0$
%the set $\Omega_{u}$ of $\sm\in\Specm \Ga$ for
%which there exists $\sn\in\Specm \Ga$, such that the image of $u$
%in $\cA(\sm,\sn)$ is nonzero
%
there exists a massive set $\Omega_{u}\subset\Specm \Ga$
such that for every $\sm\in \Omega_{u}$ the image $\bar{u}$
of $u$ in $U/U\sm$ is nonzero.
\end{proposition}

\begin{proof}
 Let $\sm\in \Specm
\Ga$, $N=u\Ga$.
%Then the exact sequence of right $\Ga$-modules
%\begin{equation}
%\label{equation-annihilating-sequence}
%0\longrightarrow N\longrightarrow U %\longrightarrow U/N \longrightarrow0
%\end{equation}
%becomes non-exact after the tensoring with %$\Ga/\sm$. We reduce this situation to a problem %about finitely generated $\Ga$-modules.
Then $\bar{u}=0$ if and only if $u=\su {i=1}{n}u_{i}m_{i},$ for some $ u_{i}\in U, m_{i}\in \sm,\ i=1,\dots,n.$
Denote $S=\cu{i=1}{n} \supp u_{i}$ and $M=U(S)$. If $\bar{u}=0$, then the
exact sequence of right $\Ga$-modules
\begin{equation*}
%\label{equation-annihilating-sequence-of-}
0\longrightarrow N\longrightarrow M \longrightarrow M/N \longrightarrow0
\end{equation*}
becomes non-exact after tensoring with $\Ga/\sm$, i.e.
$\Tor_{1}^{\Ga}(M/N,\Ga/\sm)\ne0$. Denote $$X(u,S)=\{\sm\in \Specm \Ga\mid \Tor_{1}^{\Ga}(M/N,\Ga/\sm)=0\}.$$ Following Lemma
\ref{lemma-tor-is-nonzero-on-closed-set} we can set
$\Omega_u=\ca {S\subset \cM}{} X(u,S)$.
\end{proof}

\section{Representations of    Galois orders}\label{section-repres-integral}

\subsection{Extension of characters for Galois orders}
\label{subsection-Finiteness-of-extensions-for-integral-Galois-algebras}

Given $\sm\in\Specm \Ga$ fix $\ell_{\sm}\in \Specm \bar{\Ga}$ and denote $\cM_{\ell_{\sm}}\subset\cM$ ($G_{\ell_{\sm}}\subset G$) the stabilizer of $\ell_{\sm}$ in $\cM$ (in $G$). The ideal $\sm$ defines $\cM_{\ell_{\sm}}$ and $G_{\ell_{\sm}}$ uniquely up to $G$-conjugation. It allows us to use the notation $\cM_{\sm}$  instead of  $\cM_{\ell_{\sm}}$ and ${G}_{\sm}$  instead of  ${G}_{\ell_{\sm}}$.

 Denote by
$S(\sm,\sn)$ the following $G$-invariant subset in $\cM$
\begin{equation}
\label{equation-def-of-s-m-n}
S(\sm,\sn)=\{m\in\cM\,|\,\ell_{\sn}\in GmG\cdot \ell_{\sm} \}=
\{m\in\cM\,|\,g_{2}\ell_{\sn}
\end{equation}
$$=mg_{1}\ell_{\sm}\text{ for some } g_{1},g_{2}\in G \}.
$$
Obviously this definition does not depend on the choice of
$\ell_{\sm}$ and $\ell_{\sn}$.

\begin{lemma}\label{lemma-finite-S}
Let $\sm\in \Specm \Ga$ and $\cM$ is a group. Then
\begin{enumerate}
\item\label{enumerate-s-finiteness-by-stabilizer}
    $|\cM_{\sm}|<\infty$ if and only if for
    some $\sn\in \Specm \Ga$ at least one from the sets $S(\sm,\sn),$
    $S(\sn,\sm)$ is nonempty and finite. If  $|\cM_{\sm}|<\infty$ and  $S(\sm,\sn)\ne\varnothing,$
    $S(\sn,\sm)\ne\varnothing,$ then
    $|\cM_{{\sn}}|<\infty$ and $|S(\sm,\sn)|=|S(\sn,\sm)|.$

\item\label{enumerate-estimate-S-finiteness-of-connected}
    \begin{equation*}
     |S(\sm,\sn)|\leq
    \dfrac{|G|^{2}|\cM_{\sm }|}{|G_{\sm}||G_{{\sn}}|}.
    \end{equation*}

\item\label{enumerate-estimate-S-by-G-finiteness-of-connected}
    \begin{equation*}
    |S(\sm,\sn)/G|\leq |\{m\in \cM\mid \pi(m\ell_{\sm})=\sn\}|.
    \end{equation*}

\end{enumerate}
\end{lemma}

\begin{proof}
$S(\sm,\sn)$ is infinite if and only if for some fixed $g\in G,$ $\ell_{\sm}$ and $\ell_{\sn}$ the set $$\cM(\ell_{\sm},\ell_{\sn},g)=\{m\in\cM\mid m g
\ell_{\sm}=\ell_{\sn}\} $$  is infinite. In this case for any $m'\in \cM(\ell_{\sm},\ell_{\sn},g)$ holds $$(m'g)^{-1}\cM(\ell_{\sm},\ell_{\sn},g)g\subset \cM_{{\sm}},\ \cM(\ell_{\sm},\ell_{\sn},g)m^{\prime -1}\subset\cM_{{\sn}}.$$
To prove the equality note, that $m\in S(\sm,\sn)$ if and only if $m^{-1}\in S(\sn,\sm)$. It proves
\eqref{enumerate-s-finiteness-by-stabilizer}.
If  $|\cM_{\sm}|<\infty$ then there
exists at most $|\cM_{\sm}|$ elements $\vi\in\cM$, such
that $\ell_{\sn}=\vi\ell_{\sm}$. Considering the $G$-orbits of
$\ell_{\sm}$ and $\ell_{\sn}$, which have the lengths $\dfrac{|G|}{|G_{\sm}|}$ and $\dfrac{|G|}{|G_{{\sn}}|}$ respectively, we  obtain the
 inequality  \eqref{enumerate-estimate-S-finiteness-of-connected}. Assume $m g_{1}\ell_{\sm}=g_{2}\ell_{\sn},$ $g_{1},g_{2}\in G,$ $m\in \cM.$ Then $(g_{1}^{-1}m g_{1}) \ell_{\sm}=$ $ g_{1}^{-1}g_{2}\ell_{\sn},$ which proves \eqref{enumerate-estimate-S-by-G-finiteness-of-connected}.
\end{proof}

The property of the Galois ring $U$ to be a Galois order
has the following
immediate impact on the representation theory of $U$. We will
consider right  Galois orders. The case of left orders is
analogous.

\begin{lemma}
\label{lemma-tight-presentation-of-degree-one}
Let $U$ be a Galois ring over $\Gamma,$  $\Ga$  a noetherian algebra which is a right Harish-Chandra subalgebra of $U$, $\sm\in \Specm
\Gamma,$ such that $\cM_{\sm}$ is finite, $S=S(\sm,\sm).$  If  $U=U\sm$ then for every $k\geq1$ there exist $\gamma_{k}\in \Gamma\setminus \sm$, $v_{j}\in U$, $\nu_{j}\in \sm^{k}$, $j=1,\dots,N$ such that
\begin{equation}
\label{equation-inductive-gamma-k}\gamma_{k}=\sum_{j=1}^{N} v_{j}\nu_{j}
\end{equation}
and $\supp v_{j}\in S$, $j=1,\dots,n$.
\end{lemma}

\begin{proof}
The condition $U=U\sm $ is equivalent to the condition $1\in U\sm $, i.e. holds the equality
\begin{equation}
\label{equation-decomposition-of-unit}1= \sum_{i=1}^{n} u_{i}\mu_{i}, \text{ where } u_{i}\in U, \mu_{i}\in \sm.
\end{equation}
We use induction in $k$ to prove the statement of the lemma without the condition  $\supp v_{i}\in S$, $i=1,\dots,n$. The base of induction $k=1$ follows immediately from \eqref{equation-decomposition-of-unit}.
The  induction step  is obtained by substitution of the right side of \eqref{equation-decomposition-of-unit} in \eqref{equation-inductive-gamma-k}:
 if \begin{equation}
\label{equation-crud-unit}1=\sum_{j=1}^{N} w_{j}\nu_{j}, \text{ where } w_{j}\in U, \nu_{j}\in\sm^{k}, j=1,\dots,N,
\end{equation}
then
\begin{align*}
&1=\sum_{j=1}^{N} w_{j}\nu_{j}=\sum_{j=1}^{N} w_{j}\cdot1\cdot\nu_{j}= \\
&\sum_{j=1}^{N} w_{j}\big(\sum_{i=1}^{n} u_{i}\mu_{i}\big)\nu_{j}=\sum_{j=1}^{N}    \sum_{i=1}^{n}w_{j}u_{i}(\mu_{i}\nu_{j} ).
\end{align*}
It proves the induction step, since all $\mu_{i}\nu_{j}\in \sm^{k+1}.$

 Denote      $$T=\bigcup_{i=1}^{N}\supp w_{i}\setminus S.$$
%Further we act as in the proof of Lemma \ref{lemma-S-by-m-n}.
Since $T\bigcap S(\sm,\sm)=\varnothing$, the ideals
$\ell_{\sm}^{t}$ and $\ell_{\sm}$ for every $t\in T$ belong to different
$G$-orbits. By Lemma
\ref{lemma-distinguish-orbits} there exists $f\in \Gamma$
such that $f(\ell_{\sm})\ne f(\ell_{\sm}^{t})$ for every
$t\in T$. Without loss of generality we can assume that (Lemma \ref{lemma-relation-as-in-gz-acts-zero}, \eqref{enum-f-S-mult-on-standard})
$$f_{T}(\sm,\sm)=\ds\prod_{t\in T} (f(\ell_{\sm})-
f^{t^{-1}}(\ell_{\sm}))=\ds\prod_{t\in T} (f(\ell_{\sm})-
f(\ell_{\sm}^{t}))=1.$$ In particular, it implies that $f_{T}\in 1+\sm
\otimes \Gamma +\Gamma \otimes \sm$.
 Then by Lemma \ref{lemma-relation-as-in-gz-acts-zero}, \eqref{enum-f-S-mult-on-standard}, $\gamma_{k}=f_{T}\cdot 1\in 1+\sm$. Besides, by Lemma \ref{lemma-relation-as-in-gz-acts-zero}, \eqref{enum-f-complement-to-S-throws-U-S},  $v_{i}=f_{T}\cdot w_{i}$ belongs to $U(S) $. Applying $f_{T}$ to the equality \eqref{equation-crud-unit} we obtain
$$\gamma_{k}= \sum_{i=1}^{N}v_{i}\nu_{i},
\textrm{ where }\ga_{k}\in\Ga\setminus \sm,\, v_{i}\in U(S), \nu_{i}\in \sm^{k}.$$
\end{proof}

\begin{corollary}
\label{corollary-existence-of-lifting}
%\label{enum-exists-lift-existence}
Let  $\Ga$ be a noetherian algebra,  $U$   a right  Galois order over   $\Ga$, $\sm\in \Specm
\Gamma$, such that $|\cM_{\sm}|<\infty.$ Then $U\sm\ne U$. In particular, there exists a  simple left Gelfand-Tsetlin $U$-module $M$, generated by $x\in M$, such that $\sm\cdot x=0.$
\end{corollary}

\begin{proof}
By Lemma \ref{lemma-finite-S}, \eqref{enumerate-s-finiteness-by-stabilizer} the set $S=S(\sm,\sm)$ is finite. Since  $U$ is a Galois order, by  Theorem \ref{theorem-integrality-equivalence-of-defs}, \eqref{enum-equiv-of-defs-int-right-any-s}  $U(S)$ is finitely generated as a right $\Ga$-module. Then applying the Artin-Rees Lemma (Theorem 8.5. \cite{mat})
for a right $\Ga$- module $U(S)$ and its submodule $\Ga$ we conclude that
there exists $c\geq0$ such that for every $k\geq c$
$$ U(S)\sm^{k}\cap \Gamma =(U(S)\sm^{c}\cap \Gamma )\sm^{k-c}.$$
But by Lemma~\ref{lemma-tight-presentation-of-degree-one} for every $k> c$ there exists  $\ga_{k}\in U(S)\sm^{k}\cap \Gamma$, such that $\ga_{k}\not\in\sm.$  Hence $\ga_{k}\not\in (U(S)\sm^{c}\cap \Gamma )\sm^{k-c}$, provided that $k-c>0$. The obtained contradiction shows that $U\ne U\sm$.

Denote by $v$ the image of $1$ in $U/U\sm\ne0$. Then $\sm v=0$ which
defines a structure of Gelfand-Tsetlin module on $U/U\sm$ (Lemma \ref{lemma-corollaries-of-hc-conditions}, \eqref{enumerate-induced-from-f-d-is-GZ}). We can choose by $M$ a
simple quotient of $U/U\sm$ satisfies the statement. Therefore, as the module $M$ we can choose any simple quotient of  $U/U\sm$ and set  $x$ the class of $v.$

\end{proof}

\subsection{Finiteness of extensions of characters for Galois orders}
\label{subsection-Finiteness-of-liftings-of-characters-for-Galois-orders}

 In this subsection we assume  that $U$ is a Galois
order over $\Ga$ where $\Ga$ is
normal  noetherian over $\k$. In particular, $\Gamma=\tilde{\Gamma}=U_{e}$ and
$\bar{\Ga}$ is finite over $\Ga$ by
Proposition~\ref{basics-integral}. Also $\Ga$ is a Harish-Chandra
subalgebra by Proposition~\ref{prop-integral-Galois-algebra}. If
$\ell\in \Specm \bar{\Ga}$ projects to $\sm\in \Specm \Gamma$ then
we will write $\ell=\ell_{\sm}$ and say that $\ell_{\sm}$ is lying
over $\sm$. Let $\ell_{\sm}$ and $\ell_{\sn}$ be some maximal ideals
of $\bar{\Ga}$ lying over $\sm$ and $\sn$ respectively. Note that
given $\sm\in \Specm \Ga$ the number of different $\ell_{\sm}$ is
finite due to Corollary~\ref{cor-finite-extension}.

\begin{lemma}
\label{lemma-S-by-m-n} Let  $\sm,\sn\in\Specm \Ga$,
$S=S(\sm,\sn)$, $m,n\geq 0$. Then $U=U(S)+\sn^{n}U+U\sm^{m}.$
Moreover, for every $u\in U, k\geq0$ there exists $u_{k}\in U(S)$,
such that $u\in u_{k}+\sn^{[k/2]}u\Ga+\Ga u\sm^{[k/2]}$.
\end{lemma}

\begin{proof}
Fix $u\in U$ and denote $T=\supp u\setminus S$. If
$T=\varnothing$ then $u\in U(S)$. Let $T\ne\varnothing$. We
show by induction in $k$, that there exists $u_{k}\in
U(S)$, such that
\begin{equation}
\label{equation-S-approximation}
u\in u_k+ \displaystyle
{\sum_{i=0}^k}{\sn}^{k-i} u
{\sm}^i, \ u_{k}\in U(S)\ (\text{ hence } \ u_{k+1} - u_k\in
\displaystyle {\sum_{i=0}^{k}}
{\sn}^{k-i} u {\sm}^i).
\end{equation}

Since $\ell_{\sm}^{t}$ and $\ell_{\sn}$ belong to different
$G$-orbits if $t\not\in S$, then  by Lemma
\ref{lemma-distinguish-orbits} there exists $f\in \Gamma$
such that $f(\ell_{\sn})\ne f(\ell_{\sm}^{t})$ for every
$t\in T$. Without loss of generality we can assume that
$f_{T}(\sn,\sm)=\ds\prod_{t\in T} (f(\ell_{\sn})-
f^{t^{-1}}(\ell_{\sm}))=1$, which implies $f_{T}\in 1+\sn
\otimes \Gamma +\Gamma \otimes \sm$. Set $u_{1}=f_{T} \cdot
u$. Then $u_{1}$ belongs to $u +\sn u \Gamma +\Gamma u \sm$
and, hence, $u\in u_{1} +\sn u \Gamma +\Gamma u \sm$.
Moreover, $u_{1}\in U(S)$ by Lemma \ref{lemma-relation-as-in-gz-acts-zero},
\eqref{enum-f-complement-to-S-throws-U-S}.
We prove the  induction step $k\Rightarrow k+1$. Writing in
\eqref{equation-S-approximation} the expression for $u$ in
the right hand side we obtain
\begin{align*}
&u\in u_k + {\sum_{i=0}^k} {\sn}^{k-i}
(u_k + {\sum_{j=0}^k} {\sn}^{k-j}
u{\sm}^j) {\sm}^i \subset \\ & u_k +
{\sum_{i=0}^k} {\sn}^{k-i} u_k
{\sm}^i+
 {\sum_{i=0}^{k+1}} {\sn}^{k+1-i} u
{\sm}^i,
\end{align*}
that  proofs
of the induction step, since $ u_k + {\sum_{i=0}^k}
{\sn}^{k-i} u_k {\sm}^i\subset U(S).$
\end{proof}

%\begin{remark}\label{rem-equiv-category A}
%We use notations from Lemma~\ref{lemma-corollaries-of-hc-conditions}.
%If ${\cA}({\sm},{\sn})$ is finitely generated as $\Ga$-bimodule
%then the category ${\cA}$ can be equivalently defined as follows
%(\cite{dfo:hc}, 1.4). Let $\sl,$ $\sm,$ $\sn\in \Specm \Ga$. To
%define the composition
%$$\cA(\sn,\sl)\times \cA(\sm,\sn)\myto \cA(\sm,\sl)$$
%considered pairs $a,b\in U$  and $m,n>0$, the left $\Ga$-module
%$M=L_{a,\sm^{m}}$ and the right $\Ga$-module $N=R_{b,\sl^{l}}.$ By
%the Harish-Chandra property both modules are finite dimensional
%over $\k$. As above $M$ and $N$ as $\Ga-$modules decompose
%\begin{align}
%\label{align-decomposition-M-N}
%&M=M_{0}\oplus M_{1}, \text{ where } \sn^{n}M_{0}=0, \sn M_{1}=M_{1},\\
%&N=N_{0}\oplus N_{1}, \text{ where } N_{0}\sn^{n}=0, N_{1}\sn
%=N_{1}
%\end{align}
%for some $n\geq 0$. These decompositions induce (not uniquely!)
%decompositions in $U$
%\begin{align}
%\label{align-decomposition-in-definition-of-product-in-cA}
%&a=a_{0}+a_{1}, \sn^{n}a_{0}\in U\sm^{m}, a_{1}\in \sn U+U\sm^{m}
%,\\ \nonumber &b=b_{0}+b_{1}, b_{0}\sn^{n}\in \sl^{l}U, a_{1}\in
%U\sn +\sl^{l}U.
%\end{align}
%Then the class of $b_{0}a_{0}$ in $\lr{\sl^{l}}{\sm^{m}}{U}$
%depends only on the classes of $a$ in $\lr{\sn^{n}}{\sm^{m}}{U}$
%and the class of $b$ in  $\lr{\sl^{l}}{\sn^{n}}{U}.$ Taking the
%limit as $l,m\myto \infty$ we obtain the product we need.
%
%Obviously, this definition is equivalent to the one given above.
%\end{remark}

\begin{theorem}
\label{theorem-finiteness-homs-H-Ch-homs}  For any $\sm,
\sn\in \Specm \Gamma$ such that $S=S(\sm,\sn)$ is finite, the completed
${\Gamma}_{\sn}- {\Gamma}_{\sm}$-bimodule ${\cA}({\sm},{\sn})$ (see \eqref{equation-hom-in-A}) is finitely
generated. Moreover, $\cA(\sm,\sn)$ coincides with the image of $U(S)$ in $\cA$. Besides,  $\cA(\sm,\sn)$ is finitely generated both as a right $\hat{\Ga}_{\sm}$-module and as a left $\hat{\Ga}_{\sn}$-module.
\end{theorem}

\begin{proof}
In view of Lemma \ref{lemma-S-by-m-n} and formula
\eqref{equation-hom-in-A} we have an embedding
\begin{align}
\label{align-inverse-limit} &{\cA}({\sm},{\sn}) \ \subset
\ds\lim_{\leftarrow n,m} U(S)/ \left(\left({\sn}^n U + U
{\sm}^m\right)\cap U(S)\right).\end{align}

Since $U(S)$ is a noetherian $\Gamma$-bimodule by
Theorem~\ref{theorem-integrality-equivalence-of-defs}, the
generators of $U(S)$ as a $\Gamma$-bimodule generate any $ U(S)/
(({\sn}^n U + U {\sm}^m)\cap U(S))$ as a $\Gamma$-bimodule, and
hence generate ${\cA}({\sm},{\sn})$ as a complete
${\Gamma}_{\sn}-{\Gamma}_{\sm}$-bimodule (\cite{mat}, Theorem
8.7). The statement, that  ${\cA}({\sm},{\sn})$ is finitely
generated both from the left and from the right, follows from
Theorem \ref{theorem-integrality-equivalence-of-defs},
\eqref{enum-equiv-of-defs-int-right-any-s} and from Theorem 8.7,
\cite{mat}. This completes the proof.
\end{proof}

Note that Theorem \ref{theorem-finiteness-homs-H-Ch-homs} and  Definition \ref{definition-of-category-cA} imply that the  embedding \eqref{align-inverse-limit} is an isomorphism.

The following is obvious (see Lemma \ref{lemma-Gamma-big}).
\begin{corollary}
\label{corollary-Gamma-is-big} In assuptions of
Theorem~\ref{theorem-finiteness-homs-H-Ch-homs}, $\Ga$
is big in ${\sm}$. In particular there exists only
finitely many non-isomorphic extensions of $\chi$ to simple
$U$-modules.
\end{corollary}

\subsection{Proof of Theorem \bf{A}} Corollary~\ref{corollary-existence-of-lifting} implies Theorem {\bf A},\eqref{fiber-nontrivial}.
 The condition
$|\cM_{\sm}|<\infty$ implies the finiteness of $S(\sm,\sm)$ (Lemma \ref{lemma-finite-S}, \eqref{enumerate-s-finiteness-by-stabilizer}). Consider $\chi:\Gamma\myto \k$
such that $\sm=\Ker \chi$. If $\Ga$ is not normal then
$\tilde{\Gamma}$ is a finite $\Ga$-module and $\chi$ admits
finitely many extensions to $\tilde{\Gamma}$, by
Corollary~\ref{cor-finite-extension}.
 Hence, it is enough to prove
the statement for normal $\Gamma$.
But in this case the statement follows from Corollary \ref{corollary-Gamma-is-big}, which completes
 the proof of  Theorem {\bf A}.

Combining Theorem {\bf A},\eqref{fiber-nontrivial} and Corollary~\ref{corollary-existence-of-lifting-implies-Galois-order} we obtain the following  module-theoretic characterization of
left and right Galois orders.

\begin{corollary}
\label{corollary-integrality-equivalent-to-lifting-of-characters}
Let $U$ be a Galois ring with respect to  a noetherian algebra
$\Ga$, $\cM$ a group and $m^{-1}(\Ga)\subset \bar{\Ga}$
($m(\Ga)\subset \bar{\Ga}$) for any $m\in \cM$. Then every character
$\chi:\Ga\to \k$ extends to a simple left (right) $U$-module if and
only if $U$ is right (left) Galois order.
\end{corollary}

\subsection{Bounds for dimensions and blocks of the category of Gelfand-Tsetlin modules}
\label{subsection-Blocks-of-category-of-Gelfand-Tsetlin-modules}
Denote by $r(\sm,\sn)$ the minimal number of generators of
$U(S(\sm,\sn))$ as a right $\Ga$-module. Since $\Ga$ is a Harish-Chandra subalgebra from $|S(\sm,\sm)|<\infty$ follows $r(\sm,\sn)<\infty$ and, by Theorem \ref{theorem-finiteness-homs-H-Ch-homs}, $\cA(\sm,\sm)$ is finitely generated as a right $\hat{\Ga}_{\sm}$-module (in terms of \cite{dfo:hc}, $\Ga$ is \emph{big} in $\sm$). In particular there exists finitely many non-isomorphic simple $\cA$-modules $M$ such that $M(\sm)\ne0$, besides $M(\sm)$ is finite dimensional (Lemma \ref{lemma-Gamma-big}).

\begin{lemma}
\label{lemma-ze-notion-or-rozryv}
Let
    $\sm,\sn\in\Specm \Ga$, $S=S(\sm,\sn)$, $M=U\otimes_{\Ga}\Ga/\sm$, $x=1\otimes (1+\sm)\in M(\sm)$,
    $\pi_{\sn}:M\to M(\sn)$ the canonical projection.
    Then $$\cA(\sm,\sn)\cdot x=
    \pi_{\sn}(Ux)=\pi_{\sn}(U(S)x),$$
    and
     $$\dim_{\k}M(\sn)\leq \dim_{\k}(U(S)x), \dim_{\k}M(\sn)\leq r(\sm,\sn).$$
Analogous statements hold for any $U$-module $N$
 generated by a nonzero $y\in N(\sm)$ such that $\sm y=0.$
\end{lemma}

\begin{proof}
The first equality follows from Lemma \ref{lemma-action-of-the-class}. To prove the second equality consider some $u\in U$. Then by the Chinese remainder theorem there exists $\ga\in \Ga$ such that $\ga u
x=\pi_{\sn}(ux)$ and we replace $u$ by $\ga u$. Then in the notations of Lemma \ref{lemma-S-by-m-n}, set  $k=2t$, where $\sn^{t}ux=0$. Then $\pi_{\sn}(ux)=\pi_{\sn}(u_{k}x),$ where $u_{k}\in U(S)$. The second statement is obvious.
\end{proof}

\begin{theorem}
\label{theorem-finite-dimensionity-of-component} Let $U$ be a Galois order over $\Ga$ where $\Ga$ is a
normal noetherian $\k$-algebra, $M$ is a $U$-module.
\begin{enumerate}
\item\label{enumerate-for-any-HC-module-weight-are-finite}
    If for some $\sm\in D$, $\cM_{\sm}$ is
    finite and  $M\in \bH(U,{\Gamma},D)$ is simple then $M(\sm)$ is finite
    dimensional. Both $\dim_{\k} M(\sm)$ and the number
    of isomorphism classes of simple modules $N$, such that
    $N(\sm)\ne0$, are bounded by $r(\sm,\sm)$.
\item\label{enumerate-for-any-HC-module-weight-guaranteed-by-group}
    If in addition $\cM$ is a group then for any
    $\sn\in D$  $$\dim_{\k}M(\sn)\leq
    r(\sm,\sn)<\infty.$$
\item\label{enumerate-for-free-any-root-spaces-bounded}
    Let $U$ be free as a right $\Ga$-module. If  $M$ is generated over $U$ by $x\in M$, $\sm x=0$, $\sm\in\Specm \Ga$, then
    \begin{equation*}
%    \label{equation-dimension-in-free-case}
    \dim_{\k}M(\sn) \leq |{S(\sm,\sn)}/{G}|
    \end{equation*}
\end{enumerate}
\end{theorem}

\begin{proof}
The  statements
\eqref{enumerate-for-any-HC-module-weight-are-finite} and \eqref{enumerate-for-any-HC-module-weight-guaranteed-by-group} follow from
Lemma \ref{lemma-ze-notion-or-rozryv} and  Lemma \ref{lemma-finite-S}.
 To prove
\eqref{enumerate-for-free-any-root-spaces-bounded} consider a  free right $\Ga-$module $F,$ which covers $U(S),$
$p:F\myto U(S),$ $q:F\myto U$ is the composition of $p$ with the canonical embedding, where $S=S(\sm,\sn).$  Then the image of $$q\otimes_{\Ga}\id_{K}:
F\otimes_{\Ga}K\myto U\otimes_{\Ga}K\simeq \cK, \text{ see \eqref{align-cK-obtained-by-tensoring}}$$
coincides with $KU(S)=\cK(S)\subset
\cK.$ Following Lemma \ref{lemma-dimension-with-fixed-support}  $\dim_{K} \cK(S)=|S/G|.$ From semicontinuity of the dimension of image of mapping between free modules we obtain, that for
$$q\otimes_{\Ga}\id_{\Ga/\sm}:
F\otimes_{\Ga}\Ga/\sm\myto U\otimes_{\Ga}\Ga/\sm \simeq U/U\sm $$
holds $\dim_{\k} \Im (q\otimes_{\Ga}\id_{\Ga/\sm})\leq \dim_{K} \cK(S)=|S/G|.$
\end{proof}

Let $D$ be an equivalence class in
$\Delta$, $\sm,\sn\in D$ and $ {\vi}\in S(\sm,\sn)$.
We say, that $\sm$ and $\sn$ are
\emph{connected}{} by $\vi$ if the following two conditions hold

\begin{enumerate}
\item\label{enumerate-connected-exist-non-zero}
There exist  $[a_{-}\vi^{-1}],$  $[a_{+}\vi]\in U$, such that
$a_{-}a_{+}^{\vi^{-1}}$  and $a_{-}^{\vi}a_{+}$ do not belong to
$\{g\ell_{\sm}, g\in G\}.$
\item\label{enumerate-connected-all-are-different}
The number of elements in the set
 $\{\pi(\vi^{g}(\ell_{\sm}))\mid g\in G\}$ equals $|G/H_{\vi}|$.
\end{enumerate}

Endow $D$ with  a structure of a
non-oriented graph as follows. The vertices are  the elements of
$D$. An edge between $\sm$ and $\sn$ exists if and only if there exists some $\vi\in \cM$ that connects $\sm$ and $\sn$. The following statement is a generalization of Theorem 32, \cite{dfo:gz}.

\begin{proposition}
\label{proposition-break} If $\sm$ and $\sn$ in $D$ are connected
by  $\vi\in \cM$, then $\sm\simeq\sn$  in $\cA_{D}$.
\end{proposition}

\begin{proof}

As in Lemma \ref{lemma-distinguish-orbits}, for every
$n\geq0$ there exists a function $f_{n}\in\Ga$, such that
$$\ga_{n}\equiv 1 \mod \ell_{\sn}^{n},  \ga_{n}\equiv 0\mod  (\vi^{g}\cdot
\ell_{\sm} )^{n} \textrm{ for } g\in G,  \ell_{\sn}\ne \vi^{g}\cdot \ell_{\sm}.$$
Consider the element $x_{n}=[a_{-}\vi^{-1}] \ga_{n}^{2}[a_{+}\vi]$
whose  coefficient by the element $e\in\cM$ equals
$$ \tau_{n}=\su{g\in G/H_{\vi}}{}a_{-}^{g^{-1}}a_{+}^{(\vi^{-1})^{g}} \big(\ga_{n}^{(\vi^{-1})^{g}}\big)^{2}.$$
By the construction $\tau_{n}\not\in \sm,$ and without loss of generality we may assume that $\tau_{n}\in 1+\sm.$  Consider two sequences in $\Ga,$  $y_{n}$ and $z_{n},$ $n\geq1,$   which converge in $\cA$ to  $1_{\sm}$ and $1_{\sn}$ respectively. Then the sequences
$$ g_{n}= \ga_{n}[a_{+}\vi]y_{n}, \,  f_{n}=z_{n}[a_{-}\vi^{-1}] \ga_{n}$$
converge to  $g:\sm\myto\sn$ and
$f:\sn\myto \sm$ respectively, such that $fg=1_{\sm}+\mu,$ where $\mu$ belongs to the image of $\sm$ in $\cA(\sm,\sm).$ Since $1_{\sm}+\mu$ is invertible, we can assume  $fg=1_{\sm}.$ Analogously one
constructs $f':\sn\myto \sm,g':\sm\myto\sn,$ such that $g'f'=1_{\sn}$.
\end{proof}

Immediately from Proposition \ref{proposition-break}  follows

\begin{corollary}
\label{corollary-nonbreak-is-isomorphism} If $\sm$ and
$\sn$ belong to the same connected component of the graph $D$ then
they are isomorphic in $\mathcal{A}_{D}$.
\end{corollary}

\begin{theorem}
\label{theorem-finiteness-of-length-of-universal-module}
 Let  $D$ be a class of $\Delta$-equivalence. Suppose $\cM$ is a group,  $D$ has a finite number of connected components and for some $\sm\in D$ the group $\cM_{\sm}$ is finite. Then the module $U/U\sm$ is of finite length.
\end{theorem}

\begin{proof}
By Corollary    \ref{corollary-nonbreak-is-isomorphism} the
skeleton $\mathscr{B}_{D}$ of the category
$\mathscr{A}_{D}$ contains a finite number of objects, $\Ob \cB_{D}=\{ \sn_{1},\dots,\sn_{k}\}$.  Consider $U/U\sm$ as an element in $\cA-\mo_{d}$ and denote the $\cB_{D}$-module $M=(U/U\sm)|_{\cB_{D}}.$ Then by Theorem \ref{theorem-finite-dimensionity-of-component}, \eqref{enumerate-for-any-HC-module-weight-guaranteed-by-group},  $\dim_{k} M(\sn_{i})\leq r(\sm,\sn_{i})$ for any $i=1,\dots, k$. Since the categories $\cA-\mo_{d}$ and $\cB_{D}-\mo_{d} $   are equivalent, it completes the proof.
\end{proof}

From the proof above we obtain

\begin{corollary}
\label{corollary-general-estimate-of-length}
In assumptions of Theorem \ref{theorem-finiteness-of-length-of-universal-module} the length of $U/U\sm$ and the number of simple objects    in the block $\bH(U,\Ga, D)$   does not exceed $\su{i=1}{k} r(\sm,\sn_{i})$, where $\sn_{i}$ runs the set of representatives of the connected components of $D$.
\end{corollary}

\subsection{Further properties of Gelfand-Tsetlin modules}
\label{subsection-Gelfand-Tsetlin-modules-for-Galois-orders}
 In this subsection we assume that $U$ is a Galois
order over noetherian $\Gamma$.
In the conditions of  Theorem \ref{theorem-finiteness-homs-H-Ch-homs} we are able to prove the following generalization of Corollary \ref{corollary-existence-of-lifting}.

\begin{theorem}
\label{theorem-image-of-u(s)-in-cA-is-nonzero}
Let $U$ be a Galois order over a noetherian algebra, $\sm,\sn\in\Specm \Ga$ such that  $S(\sm,\sn)$, $|\cM_{\sm}|,$ $|\cM_{\sn}|$ are finite and $S(\sm,\sn)\ne\varnothing.$ Then the image of $U(S)$ in $\cA(\sm,\sn)$ (Definition \eqref{equation-hom-in-A}) is nonzero.
\end{theorem}

\begin{proof}
Note that for a nonempty $G$-invariant $S\subset \cM$ the $\Ga$-bimodule $U(S)$ is nonzero, since $K U(S)\subset \cK$ is nonzero.
Consider $U$ as a $\Ga\otimes \Ga$-module and denote by $I$ the ideal $\sn\otimes \Ga+\Ga\otimes\sn$ in $ \La\otimes \La.$ Then the class of $u\in U$ in $\cA(\sm,\sn)$ is $0$ if for any $N\geq0$ holds $u\in I^{N}U$.
We prove the following statement: if for some $n\geq0$ holds $U(S)\subset I^{N}U$, then $U(S)\subset I^{N}U(S).$
Assume $U(S)\subset I^{N} U$, equivalently, if $v_{1},\dots,v_{k}$ are the generators of $U(S)$ as $\Ga$-bimodule, then
\begin{equation}
\label{equation-generating-annihilated-by-sm-sn-bimodule}
v_{i}=\sum_{j=1}^{k}\nu_{ij}u_{ij},\text{ for some } \nu_{ij}\in I^{N}\subset\Ga\otimes\Ga,  u_{ij}\in U, i=1,\dots,k.
\end{equation}
Set $T=\ds\bigcup_{i=1}^{n}\supp u_{i}\setminus S$. As in  the proof of Lemma \ref{lemma-S-by-m-n} construct the element $$f_{T}=1-F\in 1+\sn\otimes \Ga+\Ga\otimes\sm$$ such that for all $v_{ij}=f_{T}\cdot u_{ij}$ holds $\supp v_{ij}\subset S$.
Applying $f_{T}$ to both sides of the equality \eqref{equation-generating-annihilated-by-sm-sn-bimodule} we obtain

\begin{equation}
\label{equation-generating-annihilated-by-sm-sn-bimodule-first-iteration}
v_{i}=F\cdot v_{i}+\sum_{j=1}^{k}\nu_{ij}v_{ij},\text{ } \nu_{ij}\in I^{N},  v_{ij}\in U(S), i=1,\dots,k.
\end{equation}

Substituting the value $v_{i}$ into the right hand side of \eqref{equation-generating-annihilated-by-sm-sn-bimodule-first-iteration}  we obtain

$$ v_{i}=F^{2}\cdot v_{i}+(F+1)\sum_{j=1}^{k}\nu_{ij}v_{ij}.$$
Iterating this procedure $N-1$ times we obtain

\begin{align*}
%\label{align-generating-annihilated-by-sm-sn-bimodule-N-th-degree}
&v_{i}=F^{N}\cdot v_{i}+(F^{N-1}+\dots+F+1)\cdot\sum_{j=1}^{k}\nu_{ij}v_{ij}, \\
& \text{ for some }\nu_{ij}\in I^{N},  v_{ij}\in U(S), i=1,\dots,k,
\end{align*}
which shows $v_{i}\in I^{N}U(S)$, since $v_{i}$ itself and all $v_{ij}$ belongs to $U(S)$ and $F^{N}\in I^{N}.$
In particular, it means that
$$U(S)=\bigcap_{n=1}^{\infty}I^{N}U(S).$$

Then by the Krull Theorem (Theorem 8.9, \cite{mat}), there exists $a\in 1+I$, such that
$$ a\cdot U(S)=0.$$
Since $\Ga\otimes\Ga$ acts on $\cK$ the action of $a$ is defined on $V=U(S)K\subset \cK$ and $a\cdot V=0$ as well.
But, following Lemma~\ref{l81}, all irreducible summands of $V$ as a $K$-bimodule are of the form $V(\vi)$ for some $\vi\in\cM$, and since $\supp V=S(\sm,\sn)$, there exist coimages $\ell_{\sm}$ and $\ell_{\sn}$, such that $\ell_{\sn}=\ell_{\sm}^{\vi}$. Note, that the $K$-bimodule $V(\vi)$ is isomorphic to  $L^{H_{\vi}}$, endowed with the structure of $K$-bimodule.
 Then $a\in 1+I\subset\Ga\otimes \Ga$ and it can be written in the form

$$a=1+\sum_{i=1}^{m} \nu_{i}\otimes \alpha_{i}+\sum_{j=1}^{n} \beta_{j}\otimes \mu_{j}, \alpha_{i},\beta_{j}\in \Ga, \mu_{i}\in \sm,\nu_{j}\in\sn.$$
Write  the action of $a$ on $1\in L^{H_{\vi}}:$

\begin{align*}
&0=a\cdot 1=(1+\sum_{i=1}^{m} \nu_{i}\otimes \alpha_{i}+\sum_{j=1}^{n} \beta_{j}\otimes \mu_{j})\cdot 1=\\
&1+\sum_{i=1}^{m} \nu_{i}^{\vi}\alpha_{i}+\sum_{j=1}^{n} \beta_{j}^{\vi}\mu_{j}\in 1+\ell_{\sn},
\end{align*}
because all  the elements in the formulas above belong to $\overline{\Ga}$, $\nu_{i}\in \ell_{\sn}$ and $\mu_{j}^{\vi}\in \ell_{\sn}$, since $\sm^{\vi}\subset\ell_{\sm}^{\vi}=\ell_{\sn}.$
But $0\not\in 1+\ell_{\sn}$, which completes the proof.
\end{proof}

Note that this theorem in the case $\sm=\sn$ together with Theorem \ref{theorem-HC-mod-equiv-A-mod} gives another proof of Corollary \ref{corollary-existence-of-lifting}.

  Let $\bar{\Ga}$ be the integral
closure of $\Ga$ in $L$, $\vi\in\cM$ and $i:\Ga\to \bar{\Ga}$,
$i_{\vi}:{\vi}(\Ga)\to \bar{\Ga}$ the canonical embeddings, $\pi$
and $\pi_{\vi}$ the induced mappings of the maximal spectra,
$\sm\in\Specm \Ga$, $\pi^{-1}(\sm)=\{\ell_{1},\dots,\ell_{k}\}$.
The following lemma describes the sets $X_{a}, a\in U$ (see Lemma~\ref{lemma-corollaries-of-hc-conditions}).

\begin{lemma}
\label{lemma-X-u-for-Galois-algebras}  Let $a\in
L^{H_{\vi}}$ and $V=\Ga [a\vi]\Ga$. Then the set of simple
factors of the left $\Ga$-module $V\otimes_{\Ga}\Ga/\sm$
coincides with the set of simples of the form $\Ga/\sn$,
$\sn\in\pi({\vi}(\pi^{-1}(\sm)))=\{\pi(\ell_{1}^{\vi}),$ $\dots,$
$\pi(\ell_{k}^{\vi}) \}$. Besides, for $u\in U$ and
$\sm\in\Specm \Ga$ holds $X_{u}(\sm)=\pi(\supp u\cdot
\pi^{-1}(\sm))$.
\end{lemma}

\begin{proof}
By Remark \ref{remark-squarebracket-gamma-generated-bimodule} $V/V\sm\simeq \Ga\Ga^{\vi}/\Ga\sm^{\vi}.$  Since $\Ga^{\vi}\subset
\Ga\Ga^{\vi}$ is a finite integral  extension, holds $(\Ga\Ga^{\vi})\sm^{\vi}=\Ga\sm^{\vi}\ne\Ga\Ga^{\vi}\subset\bar{\Ga}.$
The kernels of homomorphisms $p:\bar{\Ga}\to \k$ extending
$p^{\vi}:\Ga^{\vi}\to \Ga^{\vi}/\sm^{\vi}$ form a
$G^{\vi}$- orbit
$\{\pi(\ell_{1}^{\vi}),\dots,\pi(\ell_{k}^{\vi})\}$, whose
restrictions on $\Ga$ uniquely define all characters of $\Ga\Ga^{\vi},$ extending $p^{\vi}.$
This proves the first statement.

Let $V=\Gamma u\Gamma$. Then Lemma
\ref{lemma-relation-as-in-gz-acts-zero},
\eqref{enum-f-S-mult-split-components-mono} reduces the
second statement to the case of $\Ga$-bimodule $V$
generated by elements of the form
$[a_{1}\vi],\dots,[a_{k}\vi]$. Hence the  second statement
follows from the first one.
\end{proof}

We assume that $\Gamma$ is  normal noetherian $\k$-algebra and
$\Om_{2}$ and $\Om_{r}$ are as in Section
\ref{section-Representations-of-general-position}.
For $\sm\in \Specm \Gamma$ denote by $D(\sm)$ denote
the class of $\Delta$-equivalence of $\sm,$ where
$\Delta$  is defined in \ref{subsection-Gelfand-Tsetlin-modules}.

\begin{theorem}
\begin{enumerate}
\item\label{enum-homs-in-ha-cha-cat-no-s-n0-hom} If
    $S(\sm,\sn)=\varnothing$ for some  $\sm,\sn\in
    \Specm \Ga,$ then $\cA(\sm,\sn)=0$.

\item \label{enum-homs-in-ha-cha-cat-x-a-rel-supp}
Let $\Delta'$ be the minimal equivalence, containing all $(\sm,\sn)\in \Specm \Ga\times\Specm \Ga$ such that $S(\sm,\sn)\ne\varnothing$ Then $\Delta=$ $\Delta'.$
Besides, under the assumption of Theorem \ref{theorem-finiteness-homs-H-Ch-homs} the category $\cA(D)$  does not split into a non-trivial direct sum and acts faithfully on $\bH(U,{\Gamma},D).$

\item\label{enum-homs-in-ha-cha-cat-one-dim} If
    $\sm\in\Om_{r}$, then $\cA(\sm,\sm)$ is a
    homomorphic image  of $\ {\Ga}_{\sm}$. In
    particular, there exists a unique up to isomorphism
    simple $U$-module $M$, extending the character
    $\chi:\Ga\myto \Ga/\sm$.

\item\label{enum-homs-in-ha-cha-cat-one-dim-univ-mod}
    Let $\sm\in\Om_{r}$, $D=D(\sm)$,
    $M_{\sm}=\cA_{D}/\cA_{D}\hat{\sm}$, where
    $\hat{\sm}\subset {\Gamma}_{\sm}$ is the
    completed ideal.  Then $U/U\sm$ is canonically
    isomorphic to $\bF(M_{\sm})$. In particular, if $\sm\in\Omega_r$, then $U/U\sm$ is simple.

\item\label{enum-homs-in-ha-cha-cat-gamma-completed-in-point}
    Let $\cM$ be a group, $\sm\in\Om_{r}\cap\Om_{2}$.
    Then for every $\sn \in D(\sm)$ all objects of $\cA_{D}$ are isomorphic and 
    $$\cA(\sn,\sn)\simeq \hat{\Gamma}_{\sn}.$$
\end{enumerate}
\label{theorem-homs-in-ha-cha-cat}
\end{theorem}

\begin{proof}
The statement \eqref{enum-homs-in-ha-cha-cat-no-s-n0-hom}
follows from Lemma \ref{lemma-S-by-m-n}.
By Lemma \ref{lemma-X-u-for-Galois-algebras} $\vi\in S(\sm,\sn)$ if and only if $\Ga/\sn$ is a right subfactor of $\Ga [a\vi]\Ga/\Ga [a\vi]\sm.$ It proves the first statement from \eqref{enum-homs-in-ha-cha-cat-x-a-rel-supp}. To prove the second statement note, that $U(\sm,\sn)\ne0$ if and only if $S(\sm,\sn)\ne\varnothing$ and, following Theorem \ref{theorem-image-of-u(s)-in-cA-is-nonzero} and \eqref{enum-homs-in-ha-cha-cat-no-s-n0-hom}, if and only if $\cA(\sm,\sn)\ne0.$ On other hand, if $a\in\cA(\sm,\sn),$ $a\ne0,$ then there exists $N\geq1,$ such that $a\not\in\cA(\sm,\sn)\sm^{N},$ hence $a$ acts nontrivially on $U/U\sm^{N}.$  It proves statement on $\cA(D).$
 To prove  \eqref{enum-homs-in-ha-cha-cat-one-dim} we note
that, for $\sm\in\Omega_{r}$ holds $|S(\sm,\sm)|=1$ hence by Theorem \ref{theorem-finiteness-homs-H-Ch-homs}  $\cA(\sm,\sm)$ is generated as $\hat{\Ga}_{\sm}$-bimodule by the class of $e\in U$ $\bar{e},$ which is central element in $\cA(\sm,\sm).$ On other hand, there
exists the canonical complete algebra homomorphism
$i:\hat{\Ga}_{\sm}\myto \cA(\sm,\sm)$, $i(1)=\bar{e}$,
which is surjective.
 By Theorem \ref{theorem-HC-mod-equiv-A-mod} $\cA_{D}-\mo_{d}$
is equivalent to the full subcategory
$$\bF(\cA_{D}-\mo_{d})\simeq\bH(\Ga,U,D)\subset U-\mo.$$ For $\sm\in D$ consider the functor $W_{\sm}:\bH(\Ga,U,D)\myto \k-\Mod$
$$W_{\sm}(M)=\{x\in M\,|\,\sm\cdot x=0\}.$$
This is a representable functor, namely
$$W_{\sm}\simeq \Hom_{U}(U/U\sm,\,\ - \ ).$$ For $N=\bF(N')$ we have
$$\Hom_{U}(\bF(M_{\sm}),\bF(N'))\simeq
\Hom_{\cA}(M_{\sm},N')\simeq W_{\sm}(\bF(N'))=W_{\sm}(N),$$
where all isomorphisms are functorial, i.e, $U/U\sm\simeq
\bF(M_{\sm})$. It implies
\eqref{enum-homs-in-ha-cha-cat-one-dim-univ-mod}. Consider
in $\bF(\cA_{D(\sm)}-\mo_{d})$ a $U$-module:
$$M_{\chi,n}=
U\otimes_{\Ga}\Ga/\sm^{n}\simeq
\bigoplus_{\sn\in \cO_{\sm}} \Ga/\sn^{n}.$$ Any nonzero element
from $\hat{\Ga}_{\sm}$ acts nontrivially on $M_{\chi,n}$ for any $n$.
Thus $\cA(\sm,\sm)\simeq \hat{\Ga}_{\sm}$ by the Krull intersection
theorem (\cite{mat}, Theorem 8.10, (II)).
\end{proof}

\begin{corollary}
\label{corollary-uniseriality-in-general-position} Let
$\cM$ be a group,
$D=D(\sm)\subset \Specm \Ga$ a $\Delta(U,\Ga)$-equivalence
class of a maximal ideal $\sm\in\Om_{r}\cap\Om_{2}$. Then
the category $ \bH(U,{\Gamma},D)$ is equivalent to
the category $\hat{\Gamma}_{\sm}-\mo$.
\end{corollary}

\begin{proof}
Since all  objects in $\cA_{D}$ are isomorphic by Theorem
\ref{proposition-repr-of-schat},
the categories $\cA_{D}-\mo$ and $\cA(\sm,\sm)-\mo$ are
equivalent. Note that the functors of restriction
$res:\bH(U,{\Gamma},D)\to \cA(\sm,\sm)-\mo$ and  of
induction $ind:\cA(\sm,\sm)\to \bH(U,{\Gamma},D)$ are
quasi-inverse.
\end{proof}

\begin{remark}
\label{remark-formal-series} Recall that if $\sm$ is
nonsingular point of $\Specm \Ga$, then $\hat{\Ga}_{\sm}$
is isomorphic to the algebra of formal power series in
$\gkdim \Ga$ variables.
\end{remark}

\begin{remark}
\label{remark-nonzero-in-u-zero-in-a} In the assumption of
Theorem \ref{theorem-homs-in-ha-cha-cat},
\eqref{enum-homs-in-ha-cha-cat-x-a-rel-supp} a nonzero
element $u\in S(\sm,\sn)$ can turn zero in $\cA(\sm,\sn)$.
\end{remark}

\subsection{Proof of Theorem {\bf B}}
First statement of Theorem~B follows from Theorem~\ref{proposition-repr-of-schat}.

Theorem~\ref{proposition-repr-of-schat},
Theorem~\ref{theorem-homs-in-ha-cha-cat}  and
Corollary~\ref{corollary-uniseriality-in-general-position}
immediately imply the second part of {Theorem~B}.

The third statement of Theorem~B is an analogue of the
Harish-Chandra theorem for the universal enveloping
algebras (\cite{d}). In particular it shows that the
subcategories in $U-\mo$, described in
Corollary~\ref{corollary-uniseriality-in-general-position},
contain enough modules. Suppose that conditions of Theorem
B are satisfied. Consider the massive set $\Omega_u$
constructed in
Proposition~\ref{proposition-embedding-U-in-cA} and set
$$\Omega_{u}'=\Omega_{u}\cap \Omega_{2}\cap \Omega_{r}.$$ Then for
any $\sm\in \Omega_{u}'$ the element $u$ acts nontrivially
on
 $U/U\sm$ which is simple by
Theorem~\ref{theorem-homs-in-ha-cha-cat}. This completes
the proof of Theorem~B.

\section{Gelfand-Tsetlin modules for $\gl_n$}
\label{section-Gelfand-Tsetlin-modules-for-gl-n}

Consider the general lineal Lie algebra $\gl_n$ with the standard
basis $e_{ij}, i,j=1, \ldots, n$.  Set $U=U(\gl_n)$, $U_{m}=$
$U(\gl_{m})$, $1\leq m\leq n$. Let $Z_{m}$ be the center of $U_{m}$. We identify
$\gl_{m}$ for $m\leqslant n$ with a Lie subalgebra of $\gl_n$
spanned by $\{ e_{ij}\,|\, i,j=1,\ldots,m \}$, so that we have the
following chain of inclusions
$$\gl_1\subset \gl_2\subset \ldots \subset \gl_n.$$
It induces the inclusions $U_1\subset$ $U_2\subset$ $\ldots$
$\subset U_n$ of the universal enveloping algebras.  The
\emph{Gelfand-Tsetlin} subalgebra (\cite{dfo:gz}) ${\Gamma}$ in $U$
is generated by $\{ Z_m\,|\,m=1,\ldots, n \}$, where $Z_m$ is a
polynomial algebra in $m$ variables $\{ c_{mk}\,|\,k=1,\ldots,m
\}$,
\begin{equation}\label{equ_3}
 c_{mk } \ = \ \displaystyle {\sum_{(i_1,\ldots,i_k)\in \{
1,\ldots,m \}^k}} e_{i_1 i_2}e_{i_2 i_3}\ldots e_{i_k i_1}.
\end{equation}
Hence ${\Gamma}$ is a polynomial algebra in $\displaystyle
\frac{n(n+1)}{2}$ variables $\{ c_{ij}\,|\, 1\leqslant j\leqslant
i\leqslant n \}$ (\cite{zh:cg}).
Denote by $K$ be the field of
fractions of ${\Gamma}$. Let $\cM\subset$ ${\cL}$, $\cM\simeq$ $\bZ^{\frac{n(n-1)}{2}}$ be a free abelian group generated
by $\delta^{ij}$, $1\leqslant j\leqslant i\leqslant n-1$, $(\delta^{ij})_{kl}=1$ if $i=k$, $j=l$ and $0$ otherwise.
 For
$i=1,\ldots, n$ denote by $S_i$ the $i$-th symmetric group and set
$G=S_1\times \ldots \times S_n$. The group $G$ acts on ${\cL}$:
$(s\cdot \ell)_{ij}=$ $\ell_{s_j(i)\,j}$ for $\ell=$
$(\ell_{ij})\in$ ${\cL}$  and $s=$ $(s_i)\in G$. Also the
group $\cM$ acts on ${\cL}$ as follows
$\delta^{ij}\cdot\ell=$ $\ell+\delta^{ij}$, $\delta^{ij}\in$
 $\cM$.

Let $\Lambda$ be a polynomial algebra in variables
$\{\lambda_{ij}\,|$ $1\leqslant j\leqslant i\leqslant n \}$ and
$L$ be the fraction field of $\Lambda$.  Let
$\imath:{\Gamma}{\longrightarrow}$ $\Lambda$ be an embedding such
that
$$\imath(c_{mk}) \  = \ \sum_{i=1}^m
(\lambda_{mi}+m)^k \prod_{j\ne i} ( 1 -
\frac{1}{\lambda_{mi}-\lambda_{mj}} ).$$
 The image of $\imath$ coincides with the
subalgebra of $G-$invariant polynomials  in $\Lambda$ (\cite{zh:cg})
which we identify with $\Ga$.  The
homomorphism $\imath$ can be extended to the embedding of the
fields $ K\subset L$, where $L^{G}=K$ and $G=G(L/K)$ is the Galois group.
 An action of the group $G$ by conjugations on $\cM$ induces its action
 on $L*\cM$.

 Let $e$ be the identity element of the group $\cM$.
 Consider a linear map
$\mathrm t: U\longmapsto \cK$ such that

\begin{equation} \label{equation-def-t}
\mathrm t(e_{mm})=ee_{mm}, \ \mathrm t(e_{m\,m+1})=\sum_{i=1}^m
A_{mi}^+\delta^{mi}, \ \mathrm t(e_{m+1\,m})=\sum_{i=1}^m
A_{mi}^-(\delta^{mi})^{-1},
\end{equation} where
$$
A_{mi}^{\pm}= \mp \frac {\prod_j(\lambda_{m\pm 1,j}-
\lambda_{mi})} {\prod_{j\ne i}(\lambda_{mj}-\lambda_{mi})}.
$$

The map $\mathrm t$ is an algebra homomorphism by the Harish-Chandra theorem.
Moreover, we have

\begin{proposition}(see \cite{fo-GaI},Proposition 7.2)
\label{theorem-gl-n-is-galois}
 Let $S=\Ga\setminus \{0\}$. Then
\begin{itemize}
\item $\mathrm t$ is an embedding;
\item
 $UK=KU\simeq (L*\bZ^m)^{G}$, $m=n(n-1)/2$;
 \item $U$ is a Galois order over $\Ga$.
\end{itemize}
\end{proposition}
To estimate the number of isomorphism classes of simple Gelfand-Tsetlin modules for $U(\mathrm{gl}_{n})$ we need the following statement.
\begin{theorem} (\cite{fo3})
\label{theorem-u-gl-n-is-free-over-GT}
$U(\mathrm{gl}_{n})$ is free both as a left and as a right $\Ga$-module.
\end{theorem}

Set $$Q_n=\prod_{i=1}^{n-1} i!.$$

\begin{corollary}
\label{corollary-finiteness-for-gl-n} Let
$U=U(\mathrm{gl}_{n})$, $\Gamma\subset U$ is the
Gelfand-Tsetlin subalgebra, $D$  a $\Delta$-class, $\sm\in D$. Then
\begin{enumerate}
\item\label{enumerate-the-dimension-gl-simples} For a $U$-module $M$, such that $M(\sm)\ne0$ and $M$ is generated by some $x\in M(\sm)$ (in particular for a simple module) holds $$\dim_{\k}
M(\sm)\leq Q_n.$$
\item\label{enumerate-the-number-gl-simples} The
number of  isomorphism classes of
    simple $U$-modules $N$, such that $N(\sm)\ne0$ is always nonzero and does not exceed  $Q_n.$
\end{enumerate}
\end{corollary}

\begin{proof}
Note, that a simple $M$ such that $M(\sm)\ne0,$ is generated by any nonzero vector from $M(\sm).$
Since $U$ is a free $\Ga$-module, we can apply  Theorem \ref{theorem-finite-dimensionity-of-component}, \eqref{enumerate-for-free-any-root-spaces-bounded} and obtain as a boundary $\dim_{\k} M(\sm)\leq |S(\sm,\sm)/G|.$ On  the other hand, by Lemma \ref{lemma-finite-S},\eqref{enumerate-estimate-S-by-G-finiteness-of-connected}, the right hand side here is bounded by the cardinality of the set
$$ B= \{m\in \mathbb{Z}^{n(n-1)/2})\mid \pi(m+\ell_{\sm})=\ell_{\sm}\}.$$
Equivalently, $m\in B$ if and only if the $i$th rows of $\ell_{\sm}$ and $\ell_{\sm}+m$ differ
by a permutation from $S_{i},$ $i=1,\dots,n-1.$ It gives  us at most $|S_{1}|\cdot|S_{2}|\cdot\ldots \cdot |S_{n-1}|$ possibilities for $m\in\cM$ and implies \eqref{enumerate-the-dimension-gl-simples}.
 By Lemma \ref{lemma-simple-modules-in-vertex},  the  number of  isomorphism classes of
    simple $U$-modules $N$, such that $N(\sm)\ne0$, equals  the number of isomorphism classes of simple $\cA(\sm,\sm)$-modules, and the correspondence is given by $M\leftrightarrow M(\sm).$  Therefore, if $X=U/U\sm,$ then the $\cA(\sm,\sm)$-module $X(\sm)$ covers any simple $\cA(\sm,\sm)$-module. Together with \eqref{enumerate-the-dimension-gl-simples} it proves \eqref{enumerate-the-number-gl-simples}.
\end{proof}

\begin{remark}
We believe that the bound $Q_n$ in \eqref{enumerate-the-dimension-gl-simples} can not be improved. It is known to be exact for $n=2$ and $n=3$
\cite{dfo:gz}.
\end{remark}

\begin{remark}
Applying Theorem~{B} to the case of $U(\gl_n)$ and the Gelfand-Tsetlin subalgebra $\Ga$, we note that $\Ga$ has a simple spectrum not only on finite-dimensional modules (which is well known) but also on generic Gelfand-Tsetlin modules. On the other hand, it is known not to be the case in general (see \cite{dfo:gz}).
\end{remark}

  \section{Acknowledgment}
\noindent The first author is supported in part by the CNPq
grant (processo 301743/2007-0) and by the Fapesp grant
(processo 2005/60337-2). The second author is grateful to
Fapesp for the financial support (processos 2004/02850-2 and 2006/60763-4)
and to the University of S\~ao Paulo for the hospitality
during his visits.

\end{document}